\documentclass[11pt]{article}
\usepackage[english]{babel}
\usepackage{amsmath,amsthm,amssymb}
\usepackage{mathrsfs}
\usepackage{bbm}
\usepackage{amsfonts}
\usepackage[margin=0.8in]{geometry}
\usepackage{graphicx}
\usepackage{caption}
\usepackage{subcaption}
\usepackage{wasysym}
\usepackage{enumerate}
\usepackage{algorithm2e}
\usepackage{tabularx}
\usepackage{color}

\newtheorem{thm}{Theorem}[section]

\newtheorem{prop}[thm]{Proposition}
\theoremstyle{definition}

\theoremstyle{remark}
\newtheorem{rem}[thm]{Remark}
\numberwithin{equation}{section}

\newcommand{\R}{\mathbb R}

\newcommand{\bbF}{\mathbb F}

\newcommand{\mcF}{\mathcal F}
\newcommand{\mcI}{\mathcal I}
\newcommand{\mcJ}{\mathcal J}

\newcommand{\mcT}{\mathcal T}

\newcommand{\mcU}{\mathcal U}
\newcommand{\mcS}{\mathcal S}

\newcommand{\E}{\mathbb{E}}
\newcommand{\Prob}{\mathbb{P}}

\newcommand{\Cp}{\bar p}
\newcommand{\on}{0}
\newcommand{\off}{1}
\newcommand{\vecb}{\bold{b}}
\newcommand{\esssup}{\mathop{\rm{ess}\,\sup}}
\newcommand{\ett}{\mathbbm{1}}
\newcommand{\bigSET}{\mathcal D_{\zeta}}
\newcommand{\bigSETp}{\mathcal D_{p}}
\newcommand{\cadlag}{c\`adl\`ag~}
\newcommand{\cadlagp}{c\`adl\`ag.~}

\newcommand{\ie}{\textit{i.e.\ }}
\newcommand{\eg}{\textit{e.g.\ }}
\newcommand{\etal}{\textit{et.~al.\ }}

\begin{document}

\title{A Limited-Feedback Approximation Scheme for Optimal Switching Problems with Execution Delays\footnote{This work was supported by the Swedish Research Council through the grant NT-14, 2014-3774}}%

\author{Magnus Perninge\footnote{M.\ Perninge is with the Department of Automatic Control, Lund University, Lund,
Sweden. e-mail: magnus.perninge@control.lth.se. He is a member of the LCCC Linnaeus Center and the eLLIIT Excellence Center at Lund University.}} %
\maketitle
\begin{abstract}
We consider a type of optimal switching problems with non-uniform execution delays and ramping. Such problems frequently occur in the operation of economical and engineering systems. We first provide a solution to the problem by applying a probabilistic method. The main contribution is, however, a scheme for approximating the optimal control by limiting the information in the state-feedback. In a numerical example the approximation routine gives a considerable computational performance enhancement, when compared to a conventional algorithm.
\end{abstract}

\section{Introduction}

Consider a set of $n$ production units $F:=\{1,\ldots,n\}$ where each unit can be operated at two different levels, $\{0,1\}$, representing ``off'' and ``on''. We assume that a central operator can switch production between the two operating levels in each unit. Following a switch from ``off'' to ``on'' in Unit $i$ the output will, in general, not immediately jump to installed capacity, $\Cp_i$. Rather we assume that the production ramps up during a delay period $[0,\delta_i]$, with $\delta_i>0$. We thus assume that the output of Unit $i$ following a switch from ``off'' to ``on'' is described by a Lipschitz continuous function $R_i:[0,\delta_i]\to [0,\Cp_i]$, with $R_i(0)=0$ and $R_i(\delta_i)=\Cp_i$. Turning off the unit is, on the other hand, assumed to render an immediate halt of production.

We consider the problem where the central operator wants to maximize her return over a predefined operation period $[0,T]$ (with $T<\infty$) that can represent, for example, the net profit from electricity production in $n$ production units or mineral extraction from $n$ mines. The profit depends on the operating-mode and the output from the $n$ units, but also on an observable diffusion process $(X_t: 0\leq t\leq T)$.

For $i=1,\ldots,n$ we let $0\leq\tau^i_1\leq \cdots \leq\tau^i_{N_i}< T$ represent the times that the operator intervenes on Unit $i$. We assume, without loss of generality, that all units are off at the start of the period so that intervention $\tau^i_{2j-1}$ turns operation on, while intervention $\tau^i_{2j}$ turns operation to the ``off''-mode. We define the operating-mode $(\xi_t: 0\leq t \leq T)$ of the system to be the $\mcJ:=\{0,1\}^n$--valued process representing the evolution of the operation modes for the $n$ units. The operation-mode of Unit $i$, at time $t\in[0,T]$, is then
\begin{equation*}
(\xi_t)_i:=\sum_{j=1}^{\lceil N_i/2 \rceil}\ett_{[\tau^i_{2j-1},\tau^i_{2j})}(t),
\end{equation*}
(where $\lceil a \rceil$ is the smallest integer $k$ such that $k\geq a$) and the output of the same unit is
\begin{equation*}
p_i(t):=\sum_{j=1}^{\lceil N_i/2 \rceil}\ett_{[\tau^i_{2j-1},\tau^i_{2j})}(t)R_i\left((t-\tau^i_{2j-1})\wedge \delta_i\right),
\end{equation*}
with the convention that $\tau^i_{N_i+1}=\infty$. Each intervention on Unit $i$ renders a cost $c_i^{\on}:[0,T]\to \R_+$ when turning operation from ``off'' to ``on'' and a cost $c_i^{\off}:[0,T]\to \R$ when the intervention is turning off the unit. We assume that a given operation strategy $u:=(\tau^1_1,\ldots,\tau^1_{N_1};\ldots;\tau^n_1,\ldots,\tau^n_{N_n})$ gives the total reward
\begin{equation}\label{ekv:objFUN}
J(u):=\E\left[\int_0^T \psi_{\xi_t}\left(t,X_t,p(t)\right)dt+h_{\xi_T}\left(X_T,p(T)\right)-\sum_{i=1}^n\Big\{ \sum_{j=1}^{\lceil N_i/2 \rceil}c_i^{\on}(\tau^i_{2j-1})+\sum_{j=1}^{\lfloor N_i/2 \rfloor}c_i^{\off}(\tau^i_{2j})\Big\}\right],
\end{equation}
where, for each $\vecb:=(b_1,\ldots,b_n)\in\mcJ$, $\psi_{\vecb}:[0,T]\times \R^m \times \R^n_+\to \R$ and $h_{\vecb}:\R^m \times \R^n_+\to \R$ are deterministic, locally Lipschitz continuous functions of at most polynomial growth and $\lfloor a \rfloor$ is the largest integer $k$ such that $k\leq a$.

The problem of finding a maximizer of \eqref{ekv:objFUN} is a multi-modes optimal switching problem with execution delays. The multi-modes optimal switching problem was popularized by Carmona and Ludkovski in~\cite{CarmLud}, where they suggested an application to valuation of energy tolling agreements (see also the paper by Deng and Xia~\cite{DengXia}).

A formal solution to the multi-modes optimal switching problem, without delays, was derived by Djehiche, Hamad\`ene and Popier in~\cite{BollanMSwitch1}. The authors adopted a probabilistic approach by defining a verification theorem for a family of stochastic processes that specifies sufficient conditions for optimality. They further proved existence of a family of processes that satisfies the verification theorem and showed that these processes can be used to define continuous value functions that form solutions, in the viscosity sense, to a set of variational inequalities. El-Asri and Hamad\`ene \cite{ElAsri} extended the approach to switching problems where the switching costs are functions also of the state and proved uniqueness of the viscosity solutions.

Previous work on more general impulse control problems with execution delays include the novel paper by Bar-Ilan and Sulem~\cite{BarIlan95}, where an explicit solution to an inventory problem with uniform delivery lag is found by taking the current stock plus pending orders as one of the states. Similar approaches are taken by A{\"i}d \etal in \cite{RAid2} where explicit optimal solutions of impulse control problems with uniform delivery lags are derived for a large set of different problems and by Bruder and Pham~\cite{Bruder} who propose an iterative algorithm. {\O}ksendal and Sulem~\cite{OksenImpulse} propose a solution to general impulse control problems with execution delays, by defining an operator that circumvents the delay period.

A state space augmentation approach to switching problems with non-uniform delays and ramping is taken by Perninge and S\"oder in \cite{ObactMMOR} and by Perninge in \cite{CVctopf} where application to real-time operation of power systems is considered. In these papers numerical solution algorithms are proposed by means of the regression Monte Carlo approach (see Longstaff and Schwartz~\cite{Longstaff}), that has previously been proposed to solve multi-modes switching problems by Carmona and Ludkovski \cite{CarmLud} and by A{\"i}d \etal \cite{RAid}.

Although many approaches have been proposed to give solutions, both exact and approximate, to impulse control problems with execution delays, they either consider models where delays only enter through uniform lags, or they propose methods that become intractable for systems with many production units. A computational difficulty that arises when trying to find a maximizer of \eqref{ekv:objFUN} by augmenting the state with a suitable set of ``times since last intervention'' is the curse of dimensionality which may become apparent already with a relatively low number of production units\cite{Bertsekas}.

In this paper we take a different approach by limiting the feedback used in the optimization. This seems to be a computationally efficient approximation that does not sacrifice to much accuracy by deviating from optimality. Furthermore, we extend some of the main results in~\cite{BollanMSwitch1} to problems with non-uniform execution delays and ramping.

\section{Preliminaries\label{sec:Prel}}
Throughout, we will assume that $(X_t:0\leq t\leq T)$ is an $\R^m$-valued stochastic process, living in the filtered probability space $(\Omega,\mcF,\Prob)$, defined as the strong solution to a stochastic differential equation (SDE) as follows
\begin{align*}
dX_t&=a(t,X_t)dt+\sigma(t,X_t)dW_t,\quad t\in[0,T],\\
X_0&=x_0,
\end{align*}
where $(W_t; 0\leq t \leq T)$ is an $m$-dimensional Brownian motion whose natural filtration is $(\mcF^0_t)_{ 0\leq t\leq T}$, $x_0\in\R^m$ and $a:[0,T]\times \R^m \to \R^m$ and $\sigma:[0,T]\times \R^m \to \R^{m\times m}$ are two deterministic, continuous functions that satisfy
\begin{equation*}
|a(t,x)|+|\sigma(t,x)|\leq C(1+|x|)
\end{equation*}
and
\begin{equation*}
|a(t,x)-a(t,x')|+|\sigma(t,x)-\sigma(t,x')|\leq C|x-x'|,
\end{equation*}
for some constant $C>0$. We let $\bbF:=(\mcF_t)_{0\leq t\leq T}$ denote the filtration $(\mcF^0_t)_{ 0\leq t\leq T}$ completed with all $\Prob$-null sets.

We will use the following notations throughout the paper:
\begin{itemize}
  \item We let $\mcU$ be the set of all $u:=(\tau^1_1,\ldots,\tau^1_{N_1};\ldots;\tau^n_1,\ldots,\tau^n_{N_n})$ where the $\tau^i_j$ are $\bbF$--stopping times and define $\mcU_t:=\{(\tau^1_1,\ldots,\tau^1_{N_1};\ldots;\tau^n_1,\ldots,\tau^n_{N_n})\in\mcU: \tau^i_1\geq t, \text{ for }i=1,\ldots,n\}$.
  \item It is sometimes convenient to represent a control $u=(\tau^1_1,\ldots,\tau^1_{N_1};\ldots;\tau^n_1,\ldots,\tau^n_{N_n})\in\mcU$ by a sequence of intervention times $0\leq \tau_1<\cdots<\tau_N< T$, and a sequence of corresponding interventions $\beta_1,\ldots,\beta_N$ where $\tau_1:=\min_{(i,j)} \tau^i_j$ and $\tau_j:=\min_{(i,j)} \{\tau^i_j: \tau^i_j>\tau_{j-1}\}$, for $j=2,\ldots,N$, and $\beta_j:=\xi_{\tau_j}$, for $j=1,\ldots,N$. With this notation we may write the operation-mode in the more familiar form
\begin{equation*}
\xi_t:=\beta_0\ett_{[0,\tau_{1})}(t)+\sum_{j=1}^{N}\beta_j\ett_{[\tau_{j},\tau_{j+1})}(t),
\end{equation*}
with $\beta_0:=\bold{0}$ and using the convention that $\tau_{N+1}=\infty$.
  \item For each $\vecb\in\mcJ$ we let $\delta^{\vecb}\in\R^n$ be given by $(\delta^{\vecb})_i:=b_i \delta_i$, for $i=1,\ldots,n$ and let $D_\zeta^{\vecb}:=\{z\in\R^n: 0\leq z_i\leq \delta^{\vecb}_i, \text{ for }i=1,\ldots,n\}$. We define the sets $\mcI(\vecb):=\{j\in \{1,\ldots,n\} : b_i=1\}$, $\mcJ^{-\vecb}:=\{\vecb'\in\mcJ:\vecb'\neq\vecb\}$, $\bigSET:=[0,T]\times \cup_{\vecb\in\mcJ} (D_\zeta^{\vecb}\times \{\vecb\})$ and $\bigSETp:=[0,T]\times \cup_{\vecb\in\mcJ} (D_p^{\vecb}\times \{\vecb\})$.
  \item We define $D_p$ to be the domain of the production vector. Hence, $D_p:=\{p\in\R^n: 0\leq p_i\leq \Cp_i, \text{ for }i=1,\ldots,n\}$. Furthermore, for each $\vecb\in\mcJ$ we define $D_p^\vecb:=\{p\in\R^n: 0\leq p_i\leq b_i\Cp_i, \text{ for }i=1,\ldots,n\}$.
  \item We extend the functions $R_i$ by defining $R:\R^n\to D_p$ as $(R(z))_i:=R_i(z_i^+\wedge \delta_i)$, for $i=1,\ldots,n$, with $s^+=\max(s,0)$.
  \item For each $\vecb,\vecb'\in\mcJ$, we let $c^\vecb_i:=c^{b_i}_i$ and $c^{\vecb}_{\vecb'}:=\sum_{i=1}^n \ett_{[b_i\neq b'_i]}c^{\vecb}_i$.
  \item For each $\vecb\in\mcJ$ and each $u\in\mcU$ we extend the definition of $\xi_s$ to general initial conditions by defining the \cadlag process $(\xi^{\vecb}_s: 0\leq s\leq T)$ as $\xi^{\vecb}_s:= \vecb\ett_{[0,\tau_1)}(s)+\sum_{j=1}^N \beta_j \ett_{[\tau_j,\tau_{j+1})}(s)$.
  \item We let $\mcS^2$ be the set of all progressively measurable, continuous processes $(Z_t: 0\leq t\leq T)$ such that $\E\left[\sup_{t\in[0,T]} |Z_t|^2\right]<\infty$.
  \item We say that a family of processes $((Y_t^y)_{0\leq t\leq T}: y\in \R^k)$ is continuous in the parameter $y$ if
  \begin{equation*}
  \lim_{y'\to y}\E\left[\sup_{t\in [0,T]} |Y_t^{y'}-Y_t^y|\right]\to 0,\quad \forall y\in \R^k,
  \end{equation*}
  and use the notation $\|Y^{y'}-Y^y\|_1^{\Prob}:=\E\Big[\sup_{t\in [0,T]} |Y_t^{y'}-Y_t^y|\Big]$.
\end{itemize}
Further, we assume that:
\begin{itemize}
  \item We assume that the switching costs, $c_i^{\on}:[0,T]\to\R_+$ and $c_i^{\off}:[0,T]\to\R_+$ are Lipschitz continuous functions such that $\min_{t\in[0,T]}c_i^{\on}(t) + \min_{t\in[0,T]} c_i^{\off}(t) > 0$, for $i=1,\ldots,n$.
  \item We make the additional assumption that the terminal rewards $(h_\vecb)_{\vecb\in\mcJ}$ satisfy
  \begin{equation}\label{ekv:hconst}
  h_{\vecb}(x,p)\geq \max_{\beta\in\mcJ^{-\vecb}} \left\{ -c^\vecb_\beta(T) + h_{\beta}(x,p\wedge R(\delta^{\beta}))\right\},\quad \forall (x,p)\in \R^m\times D_p^{\vecb},
  \end{equation}
  which rules out any switching at time $T$.
\end{itemize}

To be able to consider feedback-control formulations we will, for all $t\in [0,T]$ and $x\in\R^m$, define the process $(X_s^{t,x};0\leq s\leq T)$ as the strong solution to
\begin{align*}
dX_s^{t,x}&=a(s,X^{t,x}_s)ds+\sigma(s,X^{t,x}_s)dW_s,\quad \forall s\in[t,T],\\
X^{t,x}_s&=x,\quad \forall s\in[0,t].
\end{align*}
A standard result (see \eg Theorem 6.16, p.~49 in \cite{XYZbook}) is that, for any $\theta\geq 1$, there exist constants $C^X_1>0$ and $C^X_2>0$ such that
\begin{equation}\label{ekv:Xbound1}
\E\left[\sup_{s\in [0,T]}|X^{t,x}_s|^\theta\right] \leq C^X_1 (1+|x|^{\theta})
\end{equation}
and for all $t'\in [0,T]$ and all $x'\in\R$
\begin{align}
&\E\left[\sup_{s\in [0,T]}|X^{t,x}_s-X^{t',x'}_s|^\theta\right]\leq C^X_2 (1+|x|^\theta)(|x-x'|^\theta+|t'-t|^{\theta/2}).\label{ekv:Xbound2}
\end{align}

As mentioned above we will assume that $\psi_{\vecb}$ and $h_{\vecb}$ are locally Lipschitz continuous and of polynomial growth, for all $\vecb\in \mcJ$. Hence, there exist constants $C^\psi>0$, $C^h>0$ and $\gamma\geq 1$ such that $|\psi_{\vecb}(t,x,p)|\leq C^\psi(1+|x|^\gamma)$ and $|h_{\vecb}(x,p)|\leq C^h(1+|x|^\gamma)$, for all $(x,t,p,\vecb)\in \R^m \times\bigSETp$.

Now, \eqref{ekv:Xbound1} implies that, for each $\theta\geq 1$, there are constants $C^{\psi}_1$ $(=C^{\psi}_1(\theta))$ and $C^h_1$ $(=C^{h}_1(\theta))$ such that, for all $(x,t,p,\vecb)\in \R^m\times \bigSETp$,
\begin{align}
\E\left[\sup_{s\in [t,T]} |\psi_{\vecb}(s,X_s^{t,x},p)|^\theta\right]&\leq C^{\psi}_1(1+|x|^{\gamma\theta})\label{ekv:psiBND1}
\end{align}
and
\begin{align}
\E\left[|h_{\vecb}(X_T^{t,x},p)|^\theta\right]&\leq C^h_1(1+|x|^{\gamma\theta})\label{ekv:hBND1}.
\end{align}
Hence, we have
\begin{align}
\E\left[\int_{0}^{T} \max_{\vecb\in\mcJ} \sup_{p\in D_p^{\vecb}}|\psi_{\vecb}(s,X_s^{t,x},p)|^\theta ds\right] \leq T C^{\psi}_1(1+|x|^{\gamma\theta})\label{ekv:psiBND3}
\end{align}
and in particular
\begin{align}
\E\left[\int_{t}^{T}\max_{\vecb\in\mcJ} \sup_{p\in D_p^{\vecb}}|\psi_{\vecb}(s,X_s,p)|^\theta ds\right] \leq T C^{\psi}_1(1+|x_0|^{\gamma\theta}).\label{ekv:Y0inS2}
\end{align}
Local Lipschitz continuity implies that, for every $\rho>0$, there exist $C^{\psi}_\rho,C^{h}_\rho>0$ such that,
\begin{equation*}
|\psi_{\vecb}(t,x,p)-\psi_{\vecb}(t,x',p')|^2\ett_{[|x|\vee |x'| \leq \rho]}\leq C^{\psi}_\rho (|x-x'|+|p-p'|)
\end{equation*}
and
\begin{equation*}
|h_{\vecb}(x,p)-h_{\vecb}(x',p')|^2\ett_{[|x|\vee |x'| \leq \rho]}\leq C^{h}_\rho (|x-x'|+|p-p'|)
\end{equation*}
for all $(x,t,p,\vecb)\in \R^m\times \bigSETp$ and $(x',t',p',\vecb)\in \R^m\times \bigSETp$. We thus have
\begin{align*}
\E\bigg[\sup_{s\in [0,T]} \!|\psi_{\vecb}(s,X_s^{t,x},p)-\psi_{\vecb}(s,X_s^{t',x'},p')|^2\bigg] &= \E\bigg[\sup_{s\in [0,T]}\! |\psi_{\vecb}(s,X_s^{t,x},p)-\psi_{\vecb}(s,X_s^{t',x'},p')|^2\ett_{[|X_s^{t,x}|\vee |X_s^{t',x'}|\leq \rho]}
\\
&\qquad+|\psi_{\vecb}(s,X_s^{t,x},p)-\psi_{\vecb}(s,X_s^{t',x'},p')|^2\ett_{[|X_s^{t,x}|\vee |X_s^{t',x'}|> \rho]}\bigg]
\\
&\leq \E\bigg[\sup_{s\in [t\vee t',T]} \Big\{C^{\psi}_\rho (|X_s^{t,x}-X_s^{t',x'}|+|p-p'|)
\\
&\qquad+C^\psi(2+|X_s^{t,x}|^{2\gamma}+|X_s^{t',x'}|^{2\gamma})\ett_{[|X_s^{t,x}|\vee |X_s^{t',x'}|> \rho]}\Big\}\bigg]
\\
&\leq C^{\psi}_\rho (C^X_2 (1+|x|)(|x-x'|+|t'-t|^{1/2})+|p-p'|)
\\
&\qquad + C^\psi(2+C^X_1(2+|x|^{2\gamma}+|x'|^{2\gamma}))\frac{C^X_1(2+|x|^{2\gamma}+|x'|^{2\gamma})}{\rho},
\end{align*}
where we have used Markov's inequality (see \eg Gut~\cite[p. 120]{Gut2}) in the last step. Now, since $\rho>0$ was arbitrary we get
\begin{align}
\lim_{(t',x',p')\to(t,x,p)}\E\bigg[\sup_{s\in [0,T]} |\psi(s,X_s^{t,x},p)-\psi(s,X_s^{t',x'},p')|^2\bigg]=0\label{ekv:psiBND2}.
\end{align}
and by a similar argument we have
\begin{align}
\lim_{(t',x',p')\to(t,x,p)}\E\left[|h(X_T^{t,x},p)-h(X_T^{t',x'},p')|^2\right]=0\label{ekv:hBND2}.
\end{align}

Furthermore, the Lipschitz continuity of $R_i$ implies that there is a constant $C^R>0$ such that $|R_i(t)-R_i(s)|\leq C^R|t-s|$, for all $(t,s)\in [0,\delta_i]^2$.\\


The above estimates will be used to provide a solution to the operators problem defined as:\\

\noindent\textbf{Problem 1.} Let $\mcU$ be the set of all $u:=(\tau^1_1,\ldots,\tau^1_{N_1};\ldots;\tau^n_1,\ldots,\tau^n_{N_n})$ where the $\tau^i_j$ are $\bbF$--stopping times. Find $u^*\in\mcU$, such that
\begin{equation}\label{ekv:OPTprob}
J(u^*)=\sup_{u\in\mcU} J(u).
\end{equation}

To facilitate the solution of Problem 1 we use the following proposition, which is a standard result for optimal switching problems with strictly positive switching costs.
\begin{prop}\label{prop:finSTRAT}
Let $\mcU^f$ be the set of finite strategies, \ie $\mcU^f:=\{u\in\mcU:\: \Prob\left[(\omega : \sum_{i=1}^n N_i(\omega)>k, \:\forall k>0)\right]=0\}$. Then,
\begin{equation}\label{ekv:optinUf}
\sup_{u\in \mcU} J(u)=\sup_{u\in \mcU^f} J(u).
\end{equation}
\end{prop}

\noindent\emph{Proof.} Assume that $u\in \mcU\setminus \mcU^f$ and let $B:=(\omega : \sum_{i=1}^n N_i(\omega)>k, \:\forall k>0)$, then $\Prob[B]>0$ and we have
\begin{align*}
J(u)\leq \E\bigg[ \int_0^T \max_{\vecb\in\mcJ} \sup_{p\in D_p^{\vecb}} |\psi_{\vecb}(s,X_s,p)|ds-\ett_B\sum_{i=1}^n{\lfloor N_i/2 \rfloor}\Big(\min_{t\in[0,T]}c_i^{\on}(t)+\min_{t\in[0,T]}c_i^{\off}(t)\Big)\bigg]=-\infty,
\end{align*}
since $\min_{t\in[0,T]}c_i^{\on}(t) + \min_{t\in[0,T]} c_i^{\off}(t) > 0$. Now, by \eqref{ekv:hBND1} and \eqref{ekv:Y0inS2} there is a constant $C>0$ such that $J(u)>-C$, for $u=\emptyset$ and \eqref{ekv:optinUf} follows.\qed\\

\section{Solution by state space augmentation}
The problem of finding a control that minimizes \eqref{ekv:objFUN} is non-Markovian in the state $(t,X_t,\xi_t)$ due to the delays, which prevents us from uniquely determining $p(t)$ from the operating mode $\xi_t$. To remove delays in impulse control problems with uniform delivery lags it was proposed in~\cite{BarIlan} to augment the state space with  the additional state, capacity of ``projects in the pipe''. With non-uniform delays and ramping this approach is not applicable. However, we can still apply a state space augmentation to remove the delays (see \eg \cite{Bertsekas}).

By adding the c\`adl\`ag, $\mcF_t$--adapted process $(\zeta_t:0\leq s\leq T)$ defined as
\begin{equation*}
(\zeta_t)_i:=\sum_{j=1}^{\lceil N_i/2 \rceil}\left((t-\tau^i_{2j-1})\wedge \delta_i\right)\ett_{[\tau^i_{2j-1},\tau^i_{2j})}(t)
\end{equation*}
we retain a Markov problem in the state $(t,X_t,\xi_t,\zeta_t)$. The output vector can now be written $p(t)=R(\zeta_t)$.

Consider a control $u\in \mcU$. For each $\vecb\in \mcJ$ and each $z\in D_\zeta^\vecb$ we define $(\zeta_s^{t,z,\vecb}:0\leq s\leq T)$ as
\begin{align*}
(\zeta_s^{t,z,\vecb})_i:=\,&\ett_{[b_i=0]}\sum_{j=1}^{\lceil N_i/2 \rceil}\left((s-\tau^i_{2j-1})\wedge \delta_i\right)\ett_{[\tau^i_{2j-1},\tau^i_{2j})}(s)
\\
&+\ett_{[b_i=1]}\Big\{\left((s-t+z_i)^+\wedge \delta_i\right)\ett_{[t,\tau^i_{1})}(s)+\sum_{j=1}^{\lfloor N_i/2 \rfloor}\left((s-\tau^i_{2j})\wedge \delta_i\right)\ett_{[\tau^i_{2j},\tau^i_{2j+1})}(s)\Big\}.
\end{align*}

\subsection{Verification theorem}
The following verification theorem is an adaptation of Theorem 1 in~\cite{BollanMSwitch1} to the case with execution delays:
\begin{thm}\label{thm:vfc}
Assume that there exists a family of processes $((Y^{t,z,\vecb}_s)_{0\leq s\leq T}: (t,z,\vecb)\in \bigSET)$ each in $\mcS^2$ such that $Y^{t,z,\vecb}_s$ is continuous in $(t,z)$ and
\begin{align}\nonumber
Y^{t,z,\vecb}_s:=\esssup_{\tau \in \mcT_{s}} \E\bigg[&\int_s^{\tau\wedge T}\psi_{\vecb}\left(r,X_r,R(z+(r-t)\vecb)\right)dr+\ett_{[\tau \geq T]}h_{\vecb}\left(X_T,R(z+(T-t)\vecb)\right)
\\
&+\ett_{[\tau < T]}\max_{\beta\in\mcJ^{-\vecb}}\left\{-c^{\vecb}_{\beta}(\tau)+Y^{\tau,(z+(\tau-t)\vecb)^+\wedge \delta^{\beta},\beta}_\tau\right\}\Big| \mcF_s\bigg].\label{ekv:Ydef}
\end{align}
Then $((Y^{t,z,\vecb}_s)_{0\leq s\leq T}: (t,z,\vecb)\in \bigSET)$ is unique and
\begin{enumerate}[(i)]
  \item Satisfies $Y_0^{0,0,0}=\sup_{u\in \mcU} J(u)$.
  \item Defines the sequence $(\tau_1^*,\ldots,\tau_{N^*}^*;\beta_1^*,\ldots,\beta^*_{N^*})$, where $(\tau_j^*)_{1\leq j\leq {N^*}}$ is a sequence of $\bbF$-stopping times given by
  \begin{equation}\label{ekv:tau1DEF}
  \tau^*_1:=\inf \Big\{s\geq 0:\: Y_s^{0,0,0}=\max_{\beta\in \mcJ^{-0}}\left\{-c^{0}_{\beta}(s)+Y^{s,0,\beta}_s\right\}\Big\}
  \end{equation}
  and
  \begin{equation}\label{ekv:taujDEF}
  \tau^*_j:=\inf \Big\{s \geq \tau^*_{j-1}:\: Y_s^{\tau^*_{j-1},z^*_{j-1},\beta^*_{j-1}}=\max_{\beta\in \mcJ^{-\beta^*_{j-1}}}\Big\{-c^{\beta^*_{j-1}}_{\beta}(s)+Y^{s,(z^*_{j-1}+(s-\tau^*_{j-1})\beta^*_{j-1})\wedge \delta^{\beta},\beta}_s\Big\}\Big\},
  \end{equation}
  for $j\geq 2$, and $(\beta_j^*)_{1\leq j\leq {N^*}}$ is defined as a measurable selection of
  \begin{equation*}
  \beta^*_j\in\mathop{\arg\max}_{\beta\in\mcJ^{-\beta_{j-1}^*}}\Big\{-c^{\beta_{j-1}^*}_{\beta}(\tau_j^*) +Y^{\tau^*_j,(z^*_{j-1}+(\tau^*_j-\tau^*_{j-1})\beta^*_{j-1})\wedge \delta^{\beta},\beta}_{\tau^*_j}\Big\},
  \end{equation*}
  where $z^*_{j}:=(z^*_{j-1}+(\tau^*_j-\tau^*_{j-1})\beta^*_{j-1})\wedge \delta^{\beta^*_{j}}$, with $z^*_0:=0$ and $\beta^*_0=0$; and $N^*:=\max\{j:\tau_j^*<T\}$. Then $u^*=(\tau_1^*,\ldots,\tau_{N^*}^*;\beta_1^*,\ldots,\beta_{N^*}^*)$ is an optimal strategy for Problem 1.
\end{enumerate}
\end{thm}

\noindent\emph{Proof.} Note that the proof amounts to showing that for all $(t,z,\vecb)\in \bigSET$, we have
\begin{align*}
Y^{t,z,\vecb}_s:=\esssup_{u\in \mcU_s^{f}} \E\bigg[\int_s^{T}\psi_{\xi_r^{\vecb}}\left(r,X_r,R(\zeta^{t,z,\vecb}_r)\right)dr+h_{\xi_T^{\vecb}}\left(X_T,R(\zeta^{t,z,\vecb}_T)\right)
-\sum_{j=1}^N c^{\beta_{j-1}}_{\beta_{j}}(\tau_j)\Big|\mcF_s\bigg],
\end{align*}
for all $s\in [t,T]$, where $\mcU^{f}_t$ is the subset of $\mcU^f$ with $\tau_1\geq t$, $\Prob$--a.s.~and $\beta_0=\vecb$. Then uniqueness is immediate, $(i)$ follows from Proposition \ref{prop:finSTRAT} and $(ii)$ follows from repeated use of the definition of the Snell envelope (see \eg Appendix D of Karatzas and Shreve \cite{KarShreve2} or Proposition 2 of Djehiche, Hamad\`ene and Popier \cite{BollanMSwitch1}).\\

First, define
\begin{align*}
Z_s:=Y_s^{0,0,0}+\int_0^s\psi_0(r,X_r,0)dr.
\end{align*}
Then by Proposition 2 of \cite{BollanMSwitch1} $Z_s$ is the smallest supermartingale that dominates
\begin{equation*}
\Big(\int_0^{s}\psi_0\left(r,X_r,0\right)dr+\ett_{[s=T]}h_0\left(X_T,0\right) + \ett_{[s < T]}\max_{\beta\in \mcJ^{-0}}\left\{-c^{\on}_{\beta}(s)+Y^{s,0,\beta}_s\right\}:0\leq s\leq T\Big)
\end{equation*}
and
\begin{align*}
Y_0^{0,0,0}&=\esssup_{\tau \in \mcT_0} \E\bigg[\int_0^{\tau\wedge T}\psi_0\left(r,X_r,0\right)dr+\ett_{[\tau\geq T]}h_0\left(X_T,0\right) + \ett_{[\tau < T]}\max_{\beta\in \mcJ^{-0}}\left\{-c^{0}_{\beta}(\tau)+Y^{\tau,0,\beta}_\tau\right\}\bigg]
\\
&=\E\bigg[\int_0^{\tau^*_1\wedge T}\psi_0\left(r,X_r,0\right)dr+\ett_{[\tau^*_1\geq T]}h_0\left(X_T,0\right) + \ett_{[\tau^*_1< T]}\max_{\beta\in \mcJ^{-0}}\left\{-c^{\on}_{\beta}(\tau^*_1)+Y^{\tau^*_1,0,\beta}_{\tau^*_1}\right\}\bigg]
\\
&=\E\bigg[\int_0^{\tau^*_1\wedge T}\psi_0\left(r,X_r,0\right)dr+\ett_{[\tau^*_1\geq T]}h_0\left(X_T,0\right) + \ett_{[\tau^*_1 < T]}\left\{-c^{\on}_{\beta^*_1}(\tau^*_1)+Y^{\tau^*_1,z_1^*,\beta^*_1}_{\tau^*_1}\right\}\bigg]
\end{align*}
Now suppose that, for some $j'>0$ we have, for all $j\leq j'$,
\begin{align*}
&Y_s^{\tau^*_{j-1},z^*_{j-1},\beta^*_{j-1}}=\E\bigg[\int_s^{\tau^*_j\wedge T}\psi_{\beta^*_{j-1}}\left(r,X_r,R(z^*_{j-1}+(r-\tau^*_{j-1})\beta^*_{j-1})\right)dr
\\
&+\ett_{[\tau^*_j \geq T]}h_{\beta^*_{j-1}}\left(X_T,R(z^*_{j-1}+(T-\tau^*_{j-1})\beta^*_{j-1})\right) + \ett_{[\tau^*_j < T]}\left\{-c^{\beta^*_{j-1}}_{\beta^*_{j}}(\tau^*_j)+Y^{\tau^*_j,z^*_{j},\beta^*_{j}}_{\tau^*_j}\right\}\Big| \mcF_s\bigg]
\end{align*}
under $\|\cdot\|_1^{\Prob}$, for each $\tau^*_{j-1}\leq s\leq T$. By the definition of $Y^{t,z,\vecb}_s$ in \eqref{ekv:Ydef} we have that for each $(t,z,\vecb)\in \bigSET$,
\begin{align*}
Z^{t,z,\vecb}:=\Big(Y_s^{t,z,\vecb}+\int_0^s\psi_{\vecb}\left(r,X_r,R(z+(r-t)\vecb)\right)dr:\: 0\leq s\leq T\Big)
\end{align*}
is the smallest supermartingale that dominates the process
\begin{align*}
\Big(\int_0^s\psi_{\vecb}\left(r,X_r,R(z+(r-t)\vecb)\right)dr+\ett_{[s=T]}h_{\vecb}\left(X_T,R(z+(T-t)\vecb)\right)
\\
+ \ett_{[s < T]}\max_{\beta\in\mcJ^{-\vecb}}\left\{-c^{\vecb}_{\beta}(s)+Y^{s,(z+(s-t)\vecb)^+\wedge \delta^{\beta},\beta}_s\right\}:\: 0\leq s\leq T\Big).
\end{align*}
For all $M\geq 1$, let $(G_l^M)_{1\leq l \leq M}$ be an $\epsilon(M)$--partition of $D_\zeta^\vecb$ (with $\epsilon(M)\to 0$ as $M\to\infty$) and let $(z_{l}^M)_{1\leq j\leq M}$ be a sequence of points such that $z_{l}^M\in G_l^M$, for $l=1,\ldots,M$. For $M,N\geq 1$ and $s\geq \tau^*_{j'}$, define
\begin{equation*}
\hat Y^{M,N}_s:=\sum_{\vecb\in \mcJ}\ett_{[\beta^*_{j'}=\vecb]}\sum_{k=0}^{N-1}\ett_{[kT/N\leq \tau^*_{j'} <(k+1)T/N]}\sum_{l=1}^M\ett_{[z^*_{j'}\in G_l^M]} Y_s^{kT/N,z_{l}^M,\vecb}.
\end{equation*}
Now, $\ett_{[\beta^*_{j'}=\vecb]}\ett_{[kT/N\leq \tau^*_{j'} <(k+1)T/N]}\ett_{[z^*_{j'}\in G_l^M]}\Big(Y^{kT/N,z_{l}^M,\vecb}_s+\int_{\tau^*_{j'}}^s\psi_{\vecb}(r,X_r,R(z_{l}^M+(r-kT/N)\vecb))dr\Big)$
is the product of a $\mcF_{\tau_{j'}^*}$--measurable positive r.v., $\ett_{[\beta^*_{j'}=\vecb]}\ett_{[kT/N\leq \tau^*_{j'} <(k+1)T/N]}\ett_{[z^*_{j'}\in G_l^M]}$, and a supermartingale, thus, it is a supermartingale for $s\geq \tau^*_{j'}$. Hence, as
\begin{align*}
\Big(\hat Y^{M,N}_s+\sum_{k=0}^{N-1}\ett_{[kT/N\leq \tau^*_{j'} <(k+1)T/N]}\sum_{l=1}^M \ett_{[z^*_{j'}\in G_l^M]}
\int_{\tau^*_{j'}}^s\psi_{\beta^*_{j'}}(r,X_r,R(z_{l}^M+(r-kT/N)\beta^*_{j'}))dr:\tau^*_{j'}\leq s\leq T\Big)
\end{align*}
is the sum of a finite number of supermartingales it is also a supermartingale.

By the continuity of $Y^{t,z,\vecb}_s$ in $(t,z)$ and the continuity of $R$ and $\psi_\vecb$ we get
\begin{align*}
&Y_s^{\tau^*_{j'},z^*_{j'},\beta^*_{j'}}+\int_{\tau^*_{j'}}^s\psi_{\beta^*_{j'}}(r,X_r,R(z^*_{j'} + (r-\tau^*_{j'})\beta^*_{j'}))dr=\mathop{\lim\,\inf}_{N,M\to\infty}\Big\{\hat Y^{M,N}_s
\\
&+\sum_{k=0}^{N-1}\ett_{[kT/N\leq \tau^*_{j'} <(k+1)T/N]}\sum_{l=1}^M \ett_{[z^*_{j'}\in G_l^M]} \int_{\tau^*_{j'}}^s\psi_{\beta^*_{j'}}(r,X_r,R(z_{l}^M+(r-kT/N)\beta^*_{j'}))dr\Big\}
\end{align*}
under $\|\cdot\|_1^{\Prob}$, for all $s\in[\tau^*_{j'},T]$. For all $\tau^*_{j'}\leq t\leq s $ we have
\begin{align*}
&\mathop{\lim\,\inf}_{N,M\to\infty}\Big\{\hat Y^{M,N}_t+\sum_{k=0}^{N-1}\ett_{[kT/N\leq \tau^*_{j'} <(k+1)T/N]}\sum_{l=1}^M \ett_{[z^*_{j'}\in G_l^M]}\int_{\tau^*_{j'}}^t\psi_{\beta^*_{j'}}(r,X_r,R(z_{l}^M+(r-kT/N)\beta^*_{j'}))dr\Big\}
\\
&\geq \mathop{\lim\,\inf}_{N,M\to\infty}\E\bigg[\hat Y^{M,N}_s+\sum_{k=0}^{N-1}\ett_{[kT/N\leq \tau^*_{j'} <(k+1)T/N]}\sum_{l=1}^M \ett_{[z^*_{j'}\in G_l^M]}
\\
&\qquad\qquad\qquad\cdot\int_{\tau^*_{j'}}^s\psi_{\beta^*_{j'}}(r,X_r,R(z_{l}^M+(r-kT/N)\beta^*_{j'}))dr\Big|\mcF_t\bigg]
\\
&\geq \E\bigg[\mathop{\lim\,\inf}_{N,M\to\infty}\hat Y^{M,N}_s+\sum_{k=0}^{N-1}\ett_{[kT/N\leq \tau^*_{j'} <(k+1)T/N]}\sum_{l=1}^M \ett_{[z^*_{j'}\in G_l^M]}
\\
&\qquad\qquad\qquad\cdot\int_{\tau^*_{j'}}^s\psi_{\beta^*_{j'}}(r,X_r,R(z_{l}^M+(r-kT/N)\beta^*_{j'}))dr\Big|\mcF_t\bigg],
\end{align*}
where the first part follows from the supermartingale property and the second inequality follows from Fatou's lemma. Hence, $\Big(Y_s^{\tau^*_{j'},z^*_{j'},\beta^*_{j'}}+\int_{\tau^*_{j'}}^s\psi_{\beta^*_{j'}}(r,X_r,R(z^*_{j'}+(r-\tau^*_{j'})\beta^*_{j'}))dr:\: \tau^*_{j'}\leq s\leq T\Big)$ is a supermartingale that dominates
\begin{align}\nonumber
\Big(\int_{\tau^*_{j'}}^s\psi_{\beta^*_{j'}}\left(r,X_r,R(z^*_{j'}+(r-\tau^*_{j'})\beta^*_{j'})\right)dr+\ett_{[s=T]} h_{\beta^*_{j'}}\left(X_T,R(z^*_{j'}+(T-\tau^*_{j'})\beta^*_{j'})\right)
\\
+ \ett_{[s< T]}\max_{\beta\in\mcJ^{-\beta^*_{j'}}}\Big\{-c^{\beta^*_{j'}}_{\beta}(s)+Y^{s,(z^*_{j'}+(t-\tau^*_{j'})\beta^*_{j'})\wedge \delta^{\beta},\beta}_s\Big\}:\: \tau^*_{j'}\leq s\leq T\Big).\label{ekv:dominerad}
\end{align}
It remains to show that it is the smallest supermartingale with this property. Let $(\tilde Z_s:\:0\leq s\leq T)$ be a supermartingale that dominates \eqref{ekv:dominerad} for all $s\in[\tau_j,T]$. Then, for each $(t,z,\vecb)\in \bigSET$ and $s\geq t$ we have
\begin{align*}
\ett_{[\beta^*_{j'}=\vecb]}\ett_{[\tau^*_{j'}=t]}\ett_{[z^*_{j'}=z]}\tilde Z_s \geq \ett_{[\beta^*_{j'}=\vecb]}\ett_{[\tau^*_{j'}=t]}\ett_{[z^*_{j'}=z]} \Big(\int_{t}^s\psi_{\vecb}\left(r,X_r,R(z+(r-t)\vecb)\right)dr
\\
+\ett_{[s=T]}h_{\vecb}\left(X_T,R(z+(T-t)\vecb)\right) + \ett_{[s < T]}\max_{\beta\in \mcJ^{-\vecb}}\Big\{-c^{\vecb}_{\beta}(s)+Y^{s,(z+(s-t)\vecb)\wedge \delta^{\beta},\beta}_s\Big\}\Big)
\end{align*}
which by \eqref{ekv:Ydef} gives, that
\begin{equation*}
\ett_{[\beta^*_{j'}=\vecb]}\ett_{[\tau^*_{j'}=t]}\ett_{[z^*_{j'}=z]}\tilde Z_s \geq \ett_{[\beta^*_{j'}=\vecb]}\ett_{[\tau^*_{j'}=t]}\ett_{[z^*_{j'}=z]} \Big(Y^{t,z,\vecb}_s+\int_t^s \psi_{\vecb}\left(r,X_r,R(z+(r-t)\vecb)\right)dr\Big).
\end{equation*}
Since this holds for all $(t,z,\vecb)\in \bigSET$ we get
\begin{equation*}
\tilde Z_s\geq Y^{\tau^*_{j'},z^*_{j'},\beta^*_{j'}}_s+\int_{\tau^*_{j'}}^s \psi\left(r,X_r,R(z^*_{j'}+(r-\tau^*_{j'})\beta^*_{j'})\right)dr
\end{equation*}
for all $s\geq \tau^*_{j'}$. Hence, $\Big(Y^{\tau^*_{j'},z^*_{j'},\beta^*_{j'}}_s+\int_{\tau^*_{j'}}^s \psi_{\beta^*_{j'}}(r,X_r,R(z^*_{j'}+(r-\tau^*_{j'})\beta^*_{j'}))dr: \: \tau^*_{j'}\leq s\leq T\Big)$ is the Snell envelope of \eqref{ekv:dominerad} and
\begin{align*}
&Y_s^{\tau^*_{j'},z^*_{j'},\beta^*_{j'}}=\E\bigg[\int_s^{\tau^*_{j'+1}\wedge T}\psi_{\beta^*_{j'}}\left(r,X_r,R(z^*_{j'}+(r-\tau^*_{j'})\beta^*_{j'})\right)dr
\\
&+\ett_{[\tau^*_{j'+1}\geq T]}h_{\beta^*_{j'}}\left(X_T,R(z^*_{j'}+(T-\tau^*_{j'})\beta^*_{j'})\right) + \ett_{[\tau^*_{j'+1} < T]}\left\{-c^{\beta^*_{j'}}_{\beta^*_{j'+1}}(\tau^*_{j'+1})+Y^{\tau^*_{j'+1},z^*_{j'+1},\beta^*_{j'+1}}_{\tau^*_{j'+1}}\right\}\Big| \mcF_s\bigg]
\end{align*}
under $\|\cdot\|_1^{\Prob}$. By induction we get that for each $N\geq 0$
\begin{align*}
Y^{0,0,0}_0=\E\bigg[&\int_0^{\tau^*_N \wedge T}\sum_{j=0}^{N\wedge N^*} \ett_{[\tau^*_{j}\leq r < \tau^*_{j+1}]} \psi_{\beta^*_j}(r,X_r,R(z^*_j+(r-\tau^*_j)\beta^*_j))dr
\\
&+\sum_{j=0}^{N\wedge N^*} \ett_{[\tau^*_{j}< T\leq \tau^*_{j+1}]}h_{\beta^*_j}(X_T,R(z^*_j+(T-\tau^*_j)\beta^*_j))-\sum_{j=1}^{N\wedge N^*}c^{\beta^*_{j-1}}_{\beta^*_{j}}(\tau^*_j)+\ett_{[\tau^*_N < T]}Y^{\tau^*_N,z^*_{N},\beta^*_{N}}_{\tau^*_N}\bigg],
\end{align*}
where $(\tau^*_0,\beta^*_0)=(0,\bold{0})$. Letting $N\to\infty$ while assuming that $u^*\in\mcU^f$ we find that $Y^{0,0,0}_0=J(u^*)$.\\

It remains to show that the strategy $u^*$ is optimal. To do this we pick any other strategy $\hat u:=(\hat\tau_1,\ldots,\hat\tau_{\hat N};\hat\beta_1,\ldots,\hat\beta_{\hat N})\in\mcU^f$ and let $(\hat z_j)_{1\leq j\leq \hat N}$ be defined by the recursion $\hat z_j:=(\hat z_{j-1}+(\hat \tau_j-\hat \tau_{j-1})\hat\beta_{j-1})\wedge \delta^{\hat\beta_j}$. By the definition of $Y^{0,0,0}_0$ in~\eqref{ekv:Ydef} we have
\begin{align*}
Y^{0,0,0}_0 &\geq \E\bigg[\int_0^{\hat \tau_1\wedge T}\psi_0\left(r,X_r,0\right)dr+\ett_{[\hat \tau_1\geq T]}h_0\left(X_T,0\right) + \ett_{[\hat \tau_1 < T]}\max_{\beta\in \mcJ^{-0}}\left\{-c^{\on}_{\beta}(\hat\tau_1)+Y^{\hat \tau_1,0,\beta}_{\hat \tau_1}\right\}\bigg]
\\
&\geq\E\bigg[\int_0^{\hat \tau_1\wedge T}\psi_0\left(r,X_r,0\right)dr+\ett_{[\hat \tau_1\geq T]}h_0\left(X_T,0\right) + \ett_{[\hat \tau_1 < T]}\left\{-c^{\on}_{\hat\beta_1}(\hat \tau_1)+Y^{\hat \tau_1,\hat z_1,\hat \beta_1}_{\hat \tau_1}\right\}\bigg]
\end{align*}
but in the same way
\begin{align*}
Y^{\hat \tau_1,\hat z_1,\hat \beta_1}_{\hat \tau_1} &\geq\E\bigg[\int_{\hat \tau_1}^{\hat \tau_2\wedge T}\psi_{\hat \beta_1} \left(r,X_r,R(\hat z_1+(r-\hat \tau_1)\hat\beta_1)\right)dr
\\
&+\ett_{[\hat \tau_2 \geq T]}h_{\hat \beta_1}\left(X_T,R(\hat z_1+(T-\hat \tau_1)\hat\beta_1)\right) + \ett_{[\hat \tau_2 < T]} \left\{-c^{\hat\beta_1}_{\hat\beta_2}(\hat \tau_2)+Y^{\hat \tau_2,\hat z_2,\hat \beta_2}_{\hat \tau_2}\right\}\Big|\mcF_{\hat \tau_1}\bigg],
\end{align*}
$\Prob$--a.s. By repeating this argument and using the dominated convergence theorem we find that $J(u^*)\geq J(\hat u)$ which proves that $u^*$ is in fact optimal and thus belongs to $\mcU^f$.\qed\\


\subsection{Existence}
Theorem~\ref{thm:vfc} presumes existence of the families $((Y^{t,z,\vecb}_s)_{0\leq s\leq T}: (t,z,\vecb)\in \bigSET)$. To obtain a satisfactory solution to Problem 1, we thus need to establish existence. The general existence proof (see \cite{CarmLud,BollanMSwitch1}) goes by defining a sequence $((Y^{t,z,\vecb,k}_s)_{0\leq s\leq T}: (t,z,\vecb)\in \bigSET)_{k\geq 0}$ of families of processes as
\begin{align}
Y^{t,z,\vecb,0}_s:=\E\left[\int_s^T\psi_{\vecb}\left(r,X_r,R(z+(r-t))\right)dr+h_{\vecb}\left(X_T,R(z+(T-t)\vecb)\right)\Big| \mcF_s\right]\label{ekv:Y0def}
\end{align}
and
\begin{align}\nonumber
Y^{t,z,\vecb,k}_s:=&\esssup_{\tau\in\mcT_{s}}\E\bigg[\int_s^{\tau\wedge T}\psi_{\vecb}\left(r,X_r,R(z+(r-t)\vecb)\right)dr+\ett_{[\tau\geq T]}h_{\vecb}\left(X_T,R(z+(T-t)\vecb)\right)
\\
&+\ett_{[\tau < T]}\max_{\beta\in \mcJ^{-\vecb}}\left\{-c^{\vecb}_{\beta}(\tau)+Y^{\tau,(z+(\tau-t)\vecb)^+\wedge \delta^{\beta},\beta,k-1}_\tau\right\}\Big| \mcF_s\bigg]\label{ekv:Ykdef}
\end{align}
for $k\geq 1$, and then showing that this sequence converges to a family $((\tilde Y^{t,z,\vecb}_s)_{0\leq s\leq T}: (t,z,\vecb)\in \bigSET)$ of $\mcS^2$--processes that satisfy the verification theorem. First we note that by letting $\mcU^k_t:=\{(\tau_1,\ldots,\tau_{N};\beta_1,\ldots,\beta_N)\in\mcU_t:\: N\leq k\}$ and using a reasoning similar to that in the proof of Theorem~\ref{thm:vfc} it follows that
\begin{align}
Y^{t,z,\vecb,k}_s=\esssup_{u\in \mcU_s^k} \E\bigg[&\int_s^{T}\psi_{\xi_r}\left(r,X_r,R(\zeta^{t,z,\vecb}_r)\right)dr+h_{\xi_T}\left(X_T,R(\zeta^{t,z,\vecb}_T)\right)-\sum_{j=1}^N c^{\beta_{j-1}}_{\beta_j}(\tau_j)\Big|\mcF_s\bigg],\label{ekv:Ykdef2}
\end{align}
with $\beta_0=\vecb$.\\

\begin{prop}\label{prop:YkinS2}
For each $k\geq 0$ we have:
\begin{enumerate}[a)]
  \item The process $(Y^{t,z,\vecb,k}_s:\:0\leq s\leq T)$ belongs to $\mcS^2$.
  \item The family $((Y^{t,z,\vecb,k}_s)_{0\leq s\leq T}: (t,z,\vecb)\in \bigSET)_{k\geq 0}$ is continuous in $(t,z)$.
\end{enumerate}
\end{prop}
\bigskip
\noindent \emph{Proof.} Mean square integrability can be deduced by noting that~\eqref{ekv:Ykdef2} and Doob's maximal inequality implies that there is a constant $C>0$, such that,
\begin{align*}
\E\left[\sup_{s\in [0,T]}|Y^{t,z,\vecb,k}_s|^2\right] &\leq C\E\bigg[\bigg(\int_0^T \max_{\vecb\in\mcJ}\max_{p\in D^{\vecb}_p}\psi_{\vecb}(r,X_r,p)dr + \max_{\vecb\in\mcJ}\max_{p\in D^{\vecb}_p} h_{\vecb}(X_T,p)\bigg)^2 \bigg]
\\
&\leq 2C\E\left[T\int_0^T \max_{\vecb\in\mcJ}\max_{p\in D^{\vecb}_p}|\psi_{\vecb}(r,X_r,p)|^2 dr + \max_{\vecb\in\mcJ}\max_{p\in D^{\vecb}_p} |h_{\vecb}(X_T,p)|^2 \right],
\end{align*}
for $k\geq 0$, and the right hand side is bounded by \eqref{ekv:psiBND1} and \eqref{ekv:hBND1}. Now, to show that b) holds we note that, for any control $u\in\mcU$, we have
\begin{align*}
\sup_{s\in [0, T]}|R(\zeta^{t,z,\vecb}_s(u))-R(\zeta^{t',z',\vecb}_s(u))|\leq n C^R \left(|z-z'|+|t-t'|\right),
\end{align*}
$\Prob$--a.s. Hence, with $D^{2,\vecb}_p(\rho):=\{(p,p')\in D^{\vecb}_p\times D^{\vecb}_p: |p-p'|\leq n C^R\rho\}$ we get, for all $k\geq 0$,
\begin{align*}
Y^{t,z,\vecb,k}_s-Y^{t',z',\vecb,k}_s& \leq\esssup_{u\in \mcU^k} \E\bigg[\int_0^{T}\psi_{\xi_r}\left(r,X_r,R(\zeta^{t,z,\vecb}_r)\right) -\psi_{\xi_r}\left(r,X_r,R(\zeta^{t',z',\vecb}_r)\right)dr
\\
&\qquad+h_{\xi_T}\left(X_T,R(\zeta^{t,z,\vecb}_T)\right)-h_{\xi_T}\left(X_T,R(\zeta^{t',z',\vecb}_T)\right)\Big|\mcF_s\bigg]
\\
& \leq \E\bigg[\int_0^{T}\max_{\vecb\in\mcJ}\max_{(p,p')\in D^{2,\vecb}_p \left(|z-z'|+|t-t'|\right)} |\psi_{\vecb}\left(r,X_r,p\right)-\psi_{\vecb}\left(r,X_r,p'\right)|dr
\\
&\qquad+\max_{\vecb\in\mcJ}\max_{(p,p')\in D^{2,\vecb}_p \left(|z-z'|+|t-t'|\right)} |h_{\vecb}\left(X_T,p\right)-h_{\xi_T}\left(X_T,p'\right)|\Big|\mcF_s\bigg],
\end{align*}
$\Prob$--a.s. Using symmetry we find that the same inequality holds for $Y^{t',z',\vecb,k}_s-Y^{t,z,\vecb,k}_s$. Now, by Doob's maximal inequality, there is a $C>0$ such that
\begin{align*}
\E\left[\sup_{s\in[0,T]}|Y^{t',z',\vecb,k}_s-Y^{t,z,\vecb,k}_s|^2\right] &\leq C \E\bigg[\int_0^{T}\max_{\vecb\in\mcJ}\max_{(p,p')\in D^{2,\vecb}_p \left(|z-z'|+|t-t'|\right)} |\psi_{\vecb}\left(r,X_r,p\right)-\psi_{\vecb}\left(r,X_r,p'\right)|^2 dr
\\
&\qquad+\max_{\vecb\in\mcJ}\max_{(p,p')\in D^{2,\vecb}_p \left(|z-z'|+|t-t'|\right)} |h_{\vecb}\left(X_T,p\right)-h_{\xi_T}\left(X_T,p'\right)|^2\bigg]
\end{align*}
and the right hand side goes to 0 as $(t',z')\to (t,z)$ by \eqref{ekv:psiBND2} and \eqref{ekv:hBND2}.\\

It remains to show that, for each $(t,z,\vecb)\in\bigSET$, the process $(Y^{t,z,\vecb,k}_s: 0\leq s\leq T)$ is continuous for all $k\geq 0$. To do this we will also show that for each $k\geq 0$ and each $(t,z,\vecb)\in\bigSET$:
\begin{enumerate}[a)]\setcounter{enumi}{2}
  \item For each $\vecb'\in\mcJ^{-\vecb}$ the process $(Y^{s,(z+(s-t)\vecb)^+\wedge\delta^{\vecb'},\vecb',k}_s: 0\leq s\leq T)$ is continuous.
\end{enumerate}
First consider the case $k=0$. We have
\begin{align*}
Y^{t,z,\vecb,0}_s=\,&\E\left[\int_0^T\psi_{\vecb}\left(r,X_r,R(z+(r-t)\vecb)\right)dr+h_{\vecb}\left(X_T,R(z+(T-t)\vecb)\right)\Big| \mcF_s\right]
\\
&-\int_0^s\psi_{\vecb}\left(r,X_r,R(z+(r-t)\vecb)\right)dr.
\end{align*}
Hence, $(Y^{t,z,\vecb,0}_s: 0\leq s\leq T)$ is the sum of a continuous process and a martingale w.r.t.~the Brownian filtration and is thus continuous. Furthermore, for all $s\leq s'\leq T$ and all $\vecb'\in\mcJ$,
\begin{align*}
|Y^{s,(z+(s-t)\vecb)^+\wedge\delta^{\vecb'},\vecb',0}_s - Y^{s',(z+(s'-t)\vecb)^+\wedge\delta^{\vecb'},\vecb',0}_{s'}|
&\leq |Y^{s,(z+(s-t)\vecb)^+\wedge\delta^{\vecb'},\vecb',0}_s - Y^{s,(z+(s-t)\vecb)^+\wedge\delta^{\vecb'},\vecb',0}_{s'}|
\\
&\qquad+ |Y^{s,(z+(s-t)\vecb)^+\wedge\delta^{\vecb'},\vecb',0}_{s'} - Y^{s',(z+(s'-t)\vecb)^+\wedge\delta^{\vecb'},\vecb',0}_{s'}|.
\end{align*}
Hence, continuity of $(Y^{s,(z+(s-t)\vecb)^+\wedge\delta^{\vecb'},\vecb',0}_s: 0\leq s\leq T)$ follows from continuity of $(Y^{t,z,\vecb,0}_s: 0\leq s\leq T)$ and continuity of $\psi$, $h$ and $R$. Moving on we assume that a)--c) hold for some $k\geq 0$. The process $\Big(Y^{t,z,\vecb,k+1}_s + \int_0^s \psi_{\vecb}\left(r,X_r,R(z+(r-t)\vecb)\right)dr : 0\leq s\leq T\Big)$ is the Snell envelope of the process
\begin{align*}
\bigg(\int_0^s\psi_{\vecb}\left(r,X_r,R(z+(r-t)\vecb)\right)dr+\ett_{[s=T]}h_{\vecb}\left(X_T,R(z+(T-t)\vecb)\right)
\\
+ \ett_{[s < T]}\max_{\beta\in \mcJ^{-\vecb}}\left\{-c^{\vecb}_{\beta}(s)+Y^{s,(z+(s-t)\vecb)^+\wedge \delta^{\beta},\beta,k}_s\right\}:\: 0\leq s\leq T\bigg)
\end{align*}
It is well known that the Snell envelope of a process $(U_s:0\leq s\leq T)$ is continuous if $U$ only has positive jumps. Now, $\Big(\int_0^s\psi_{\vecb}\left(r,X_r,R(z+(r-t)\vecb)\right)dr:\: 0\leq s\leq T\Big)$ is continuous and, since $\Big(Y^{s,(z+(s-t)\vecb)^+\wedge \delta^{\beta},\beta,k}_s:\: 0\leq s\leq T\Big)$ was assumed continuous, for all $\beta\in\mcJ$, in c),
\begin{align*}
\bigg(\ett_{[s=T]}h_{\vecb}\left(X_T,R(z+(T-t)\vecb)\right) + \ett_{[s < T]}\max_{\beta\in \mcJ^{-\vecb}}\left\{-c^{\vecb}_{\beta}(s)+Y^{s,(z+(s-t)\vecb)^+\wedge \delta^{\beta},\beta,k}_s\right\}:\: 0\leq s\leq T\bigg)
\end{align*}
is continuous on $[0,T)$ and may have a jump at $\{T\}$. By \eqref{ekv:hconst} any possible jump at time $T$ is positive, hence, $\Big(Y^{t,z,\vecb,k+1}_s:0\leq s\leq T\Big)$ is a continuous process.

By a similar argument, since $(Y^{s,(z+(s-t)\vecb)^+\wedge\delta^{\vecb'},\vecb',k}_s + \int_0^s \psi_{\vecb'}\big(r,X_r,R((z+(s-t)\vecb)\wedge \delta^{\vecb'} + (r-s)\vecb')\big)dr: 0\leq s\leq T)$ is the Snell envelope of the process
\begin{align*}
\bigg(&\int_0^s\psi_{\vecb'}\left(r,X_r,R((z+(s-t)\vecb)\wedge \delta^{\vecb'} + (r-s)\vecb')\right)dr
\\
&+\ett_{[s=T]}h_{\vecb}\left(X_T,R((z+(T-t)\vecb)\wedge \delta^{\vecb'})\right)
\\
&+ \ett_{[s < T]}\max_{\beta\in \mcJ^{-\vecb}}\left\{-c^{\vecb'}_{\beta}(s)+Y^{s,(z+(s-t)\vecb)^+\wedge \delta^{\vecb'}\wedge \delta^{\beta},\beta,k}_s\right\}:\: 0\leq s\leq T\bigg),
\end{align*}
c) holds for $k+1$. But, then a)--c) hold for $k+1$ as well. By an induction argument the proposition now follows.\qed\\

Next we show that the limiting family, $\lim_{k\to\infty}((Y^{t,z,\vecb,k}_s)_{0\leq s\leq T}: (t,z,\vecb)\in \bigSET)$, exists and satisfies the verification theorem.
\begin{thm}\label{thm:limVT}
The limit $((\tilde Y^{t,z,\vecb}_s)_{0\leq s\leq T}: (t,z,\vecb)\in \bigSET):=\lim_{k\to\infty}((Y^{t,z,\vecb,k}_s)_{0\leq s\leq T}: (t,z,\vecb)\in \bigSET)$ exists $\Prob$--a.s.~as a pointwise limit. Furthermore, the limit family $((\tilde Y^{t,z,\vecb}_s)_{0\leq s\leq T}: (t,z,\vecb)\in \bigSET)$ satisfies the verification theorem.
\end{thm}

\noindent \emph{Proof.} We need to show that the limit family $((\tilde Y^{t,z,\vecb}_s)_{0\leq s\leq T}: (t,z,\vecb)\in \bigSET)$ exists as a member of $\mcS^2$, that it is continuous in $(t,z)$ and that it satisfies \eqref{ekv:Ydef}. This is done in four steps as follows.\\

\noindent\emph{(i) Convergence.} Since $\mcU^k_t\subset \mcU^{k+1}_t$ we have that, $\Prob$-a.s.,
\begin{align*}
Y^{t,z,\vecb,k}_s \leq Y^{t,z,\vecb,k+1}_s &\leq \E\left[\int_0^T \max_{\vecb\in\mcJ}\max_{p\in D^{\vecb}_p} |\psi_{\vecb}(r,X_r,p)|dr + \max_{\vecb\in\mcJ}\max_{p\in D^{\vecb}_p} |h_{\vecb}(X_T,p)| \Big| \mcF_s\right],
\end{align*}
where the right hand side is bounded $\Prob$--a.s.~by the estimates of Section~\ref{sec:Prel}. Hence, the sequence $((Y^{t,z,\vecb,k}_s)_{0\leq s\leq T}: (t,z,\vecb)\in \bigSET)$ is increasing and $\Prob$--a.s.~bounded, thus, it converges $\Prob$--a.s.~for all $s\in [0,T]$.\\

\noindent\emph{(ii) Limit satisfies \eqref{ekv:Ydef}.} Applying the convergence result to the right hand side of \eqref{ekv:Ykdef} and using (iv) of Proposition 2 in \cite{BollanMSwitch1} we find that
\begin{align*}
\tilde Y^{t,z,\vecb}_s:=\esssup_{\tau \in \mcT_{s}} \E\bigg[&\int_s^{\tau\wedge T}\psi_{\vecb}\left(r,X_r,R(z+(r-t)\vecb)\right)dr+\ett_{[\tau \geq T]}h_{\vecb}\left(X_T,R(z+(T-t)\vecb)\right)
\\
&+\ett_{[\tau < T]}\max_{\beta\in\mcJ^{-\vecb}}\left\{-c^{\vecb}_{\beta}(\tau)+\tilde Y^{\tau,(z+(\tau-t)\vecb)\wedge \delta^{\beta},\beta}_\tau\right\}\Big| \mcF_s\bigg]\label{ekv:Ydef2}
\end{align*}

\noindent\emph{(iii) Limit in $\mcS^2$.} Using the same reasoning as above we find that there exists a constant $C>0$, such that,
\begin{align*}
\E\left[\sup_{s\in [0,T]}|\tilde Y^{t,z,\vecb}_s|^2\right] \leq C\E\left[2T\int_0^T \max_{\vecb\in\mcJ}\max_{p\in D^{\vecb}_p}|\psi_{\vecb}(r,X_r,p)|^2 dr + 2\max_{\vecb\in\mcJ}\max_{p\in D^{\vecb}_p} |h_{\vecb}(X_T,p)|^2 \right],
\end{align*}
which is bounded by the estimates of Section~\ref{sec:Prel}. To prove continuity in $s$ we note that $\tilde Y^{t,z,\vecb}_s+\int_0^s \psi_{\vecb}(r,X_r,R((z+(s-t)\vecb)^+\wedge \delta^\vecb))dr$ is the limit of an increasing sequence of continuous supermartingales and thus c\`adl\`ag~\cite{KarShreve1}. Now, for each $\vecb\in \mcJ$ and each $(t,z,\vecb)\in \bigSET$ the processes $\Big(\int_0^s \psi_{\vecb}(r,X_r,R((z+(s-t)\vecb)^+\wedge \delta^\vecb))dr: 0\leq s\leq T\Big)$ are continuous. Hence, by the properties of the Snell envelope, if $\tilde Y^{t,z,\vecb}_s$ has a (necessarily negative) jump at $s_1\in [0,T]$, then, for some $\beta_1\in\mcJ^{-\vecb}$, $\tilde Y^{s_1,(z+(s_1-t)\vecb)^+\wedge \delta^{\beta_1},\beta_1}_s$ also has a jump at $s_1$ and $\tilde Y^{t,z,\vecb}_{s_1-}=-c^{\vecb}_{\beta_1}(s_1)+\tilde Y^{s_1,(z+(s_1-t)\vecb)^+\wedge \delta^{\beta_1},\beta_1}_{s_1-}$. But, if $\tilde Y^{s_1,(z+(s_1-t)\vecb)^+\wedge \delta^{\beta_1},\beta_1}_s$ has a (negative) jump at $s_1$, then for some ${\beta_2}\in\mcJ^{-\vecb}$, the process $\tilde Y^{s_1,(z+(s_1-t)\vecb)^+\wedge \delta^{\beta_1} \wedge \delta^{\beta_2},\beta_2}_{s}$ will have a negative jump at $s_1$ and
\begin{equation*}
\tilde Y^{s_1,(z+(s_1-t)\vecb)^+\wedge \delta^{\beta_1},\beta_1}_{s_1-}=-c^{\beta_1}_{\beta_2}(s_1)+\tilde Y^{s_1,(z+(s_1-t)\vecb)^+\wedge \delta^{\beta_1} \wedge \delta^{\beta_2},\beta_2}_{s_1-}.
\end{equation*}
Repeating this argument we get a sequence $(\beta_k)_{k\geq 0}$, with $\beta_0=\vecb$ and $\beta_k\in \mcJ^{-\beta_{k-1}}$ for $k\geq 1$, such that for any $j>k\geq 0$ we have
\begin{align*}
\tilde Y^{s_1,(z+(s_1-t)\vecb)^+\wedge \delta^{\beta_1} \wedge \ldots\wedge \delta^{\beta_k},\beta_k}_{s_1-}
=-c^{\beta_k}_{\beta_{k+1}}(s_1)-\ldots-c^{\beta_{j-1}}_{\beta_{j}}(s_1)
+\tilde Y^{s_1,(z+(s_1-t)\vecb)^+\wedge \delta^{\beta_1} \wedge \ldots\wedge \delta^{\beta_j},\beta_j}_{s_1-}.
\end{align*}
Now, since $(\bigwedge_{l=1}^k \delta^{\beta_l})_{k\geq 1}$ is a decreasing sequence that takes values in a finite set and $\mcJ$ is a finite set, there are $j>k\geq 0$ such that $\bigwedge_{l=1}^j \delta^{\beta_l}=\bigwedge_{l=1}^k \delta^{\beta_l}$ and $\beta_j = \beta_k$. But then
\begin{align*}
0=-c^{\beta_k}_{\beta_{k+1}}(s_1)-\ldots-c^{\beta_{j-1}}_{\beta_{j}}(s_1)
\end{align*}
contradicting the fact that $\min_{t\in [0,T]}c^{\on}_{i}(t)+\min_{t\in [0,T]}c^{\off}_{i}(t)>0$, for all $i\in\{1,\ldots,n\}$. Hence, $\tilde Y^{t,z,\vecb}_s$ must be continuous and thus belongs to $\mcS^2$.\\

\noindent\emph{(iv) Limit continuous in $(t,z)$.} 
By the dominated convergence theorem we have,
\begin{align*}
\lim_{(t',z')\to (t,z)}\E\bigg[\sup_{s\in [0,T]}|Y^{t',z',\vecb}_s-Y^{t,z,\vecb}_s|\bigg]&=\lim_{(t',z')\to (t,z)}\E\bigg[\sup_{s\in [0,T]}\lim_{k\to\infty}|Y^{t',z',\vecb,k}_s-Y^{t,z,\vecb,k}_s|\bigg]
\\
&=\lim_{k\to\infty}\lim_{(t',z')\to (t,z)}\E\bigg[\sup_{s\in [0,T]}|Y^{t',z',\vecb,k}_s-Y^{t,z,\vecb,k}_s|\bigg] = 0.
\end{align*}
This finishes the proof.\qed\\

We have thus far derived a verification theorem for the solution of Problem 1, and shown that there exists a (unique) family of processes satisfying the verification theorem. To finish the solution of Problem 1 we show that the families of processes in the verification theorem defines continuous value functions.

\subsection{Value function representation}

We first extend the definition of the families of processes in the verification theorem to a full state-feedback form by introducing general initial conditions as follows. For all $(r,x)\in [0,T]\times \R^m$ we let $((Y^{r,x,t,z,\vecb}_s)_{0\leq s\leq T}: (t,z,\vecb)\in \bigSET)$ be the family of processes that satisfies the verification theorem for the process $(X_s^{r,x}:{0\leq s\leq T})$ and let $((Y^{r,x,t,z,\vecb,k}_s)_{0\leq s\leq T}: (t,z,\vecb)\in \bigSET)_{k\geq 0}$ be the corresponding versions of $((Y^{t,z,\vecb,k}_s)_{0\leq s\leq T}: (t,z,\vecb)\in \bigSET)_{k\geq 0}$ defined by \eqref{ekv:Y0def} and \eqref{ekv:Ykdef} with $X$ replaced by $X^{r,x}$. The following estimates hold:

\begin{prop}\label{prop:VFbnd}
There exists $C^Y_1>0$ such that, for each $\vecb \in \mcJ$, we have
\begin{equation*}
\E\left[\sup_{s\in [0,T]}|Y^{r,x,t,z,\vecb}_s|^2\right]\leq C_1^Y(1+|x|^{2\gamma}), \quad\forall \,(r,t,x,z)\in [0,T]^2\times \R^m\times D_\zeta^\vecb.
\end{equation*}
Furthermore, $Y^{r,x,t,z,\vecb}_r$ is deterministic and
\begin{equation*}
|Y^{r,x,t,z,\vecb}_r-Y^{r',x',t',z',\vecb}_{r'}|\to 0, \quad\text{as } (r',x',t',z')\to (r,x,t,z).
\end{equation*}
\end{prop}

\bigskip

\noindent\emph{Proof.} For the first part we note that, again using Doob's maximal inequality, there exists a $C>0$ such that
\begin{align*}
\E\left[\sup_{s\in [0,T]}|Y^{r,x,t,z,\vecb}_s|^2\right]&\leq C\E\left[\int_0^T \max_{\vecb\in\mcJ}\max_{p\in D_p^{\vecb}}|\psi_{\vecb}(v,X^{r,x}_v,p)|^2 dv+\max_{\vecb\in\mcJ}\max_{p\in D_p^{\vecb}}|h_{\vecb}(X^{r,x}_T,p)|^2\right]
\\
&\leq C(C^\psi_1 T+C^h_1)(1+|x|^{2\gamma}).
\end{align*}
For the second part we pick any control $u=(\tau_1,\ldots,\tau_{N};\beta_1,\ldots,\beta_N)$ $\in \mcU_r$ and let $u'=(\tau_{l}\vee r',\tau_{l+1}\ldots,\tau_{N};$ $\beta_{l},\beta_{l+1},\ldots,\beta_N)$, where $l:=\max\{j\geq 1:\tau_j\leq r'\}\vee 1$ with $\max\{\emptyset\}$=0. Then $u'¨\in\mcU_{r'}$ and we have
\begin{align*}
\sup_{s\in [0, T]}|R(\zeta^{t,z,\vecb}_s(u))-R(\zeta^{t',z',\vecb}_s(u'))|\leq n C^R \left(|z-z'|+|r-r'|+|t-t'|\right),
\end{align*}
$\Prob$--a.s. and, by Lipschitz continuity of the $c_i^{\on}$ and $c_{i}^{\off}$, the switching costs obey
\begin{equation*}
\E\bigg[\sum_{j=1}^{N'}c^{\beta'_{j-1}}_{\beta'_j}(\tau'_j) - \sum_{j=1}^{N}c^{\beta_{j-1}}_{\beta_j}(\tau_j)\bigg]=
\E\bigg[\ett_{[N>0]}c^{\vecb}_{\beta_l}(\tau'_1) - \sum_{j=1}^{l}c^{\beta_{j-1}}_{\beta_j}(\tau_j)\bigg]\leq C^c|r-r'|
\end{equation*}
for some $C^c>0$. Hence, since $u$ was arbitrary we have
\begin{align*}
&Y^{r,x,t,z,\vecb}_r-Y^{r',x',t',z',\vecb}_{r'}
\\
&\leq  \sup_{u\in\mcU}\E\bigg[\int_r^{r\vee r'}\psi_{\xi_s}\left(s,X_s^{r,x},R(\zeta^{t,z,\vecb}_s(u))\right)ds-\int_{r'}^{r\vee r'}\psi_{\vecb}\left(s,X_s^{r',x'},R(\zeta^{t',z',\vecb}_s(u'))\right)ds
\\
&\qquad +\int_{r\vee r'}^T \left\{\psi_{\xi_s}\left(s,X_s^{r,x},R(\zeta^{t,z,\vecb}_s(u))\right)-\psi_{\xi_s}\left(s,X_s^{r',x'},R(\zeta^{t',z',\vecb}_s(u'))\right)\right\}ds
\\
&\qquad +h_{\xi_T}\left(X_T^{r,x},R(\zeta^{t,z,\vecb}_T)(u)\right)-h_{\xi_T}\left(X_T^{r',x'},R(\zeta^{t',z',\vecb}_T(u'))\right)\bigg]+C^c|r-r'|.
\end{align*}
Considering \eqref{ekv:psiBND1} we see that the first two integrals on the right hand side go to zero as $r\to r'$. By arguing as in part \emph{(iv)} of the proof of Theorem~\ref{thm:limVT} we find that the remainder goes to 0 as $(r',x',t',z')\to (r,x,t,z)$. Now, by symmetry this applies to $Y^{r',x',t',z',\vecb}_{r'}-Y^{r,x,t,z,\vecb}_r$ as well and the second inequality follows.\qed\\

Repeated use of Theorem 8.5 in \cite{ElKaroui1} shows that for $k\geq 0$, there exist functions $(v_\vecb^{k})_{\vecb\in\mcJ}$ of polynomial growth, with $v_\vecb^{k}:[0,T]\times \R^m \times [0,T]\times D_\zeta^{\vecb}\to \R$ such that
\begin{equation*}
Y_s^{r,x,t,z,\vecb,k}=v_{\vecb}^k(s,X^{r,x}_s,t,z),\quad r\leq s\leq T.
\end{equation*}
Furthermore, by Theorem 8.5 in \cite{ElKaroui1} and Proposition~\ref{prop:YkinS2} the functions $v_\vecb^{k}$ are continuous. Repeating the steps in the proof of Theorem~\ref{thm:limVT} we find that the sequences $(v_\vecb^{k})^{k\geq 0}_{\vecb\in\mcJ}$ converges pointwise to functions $v_\vecb^{k}:[0,T]\times \R^m \times [0,T]\times D_\zeta^{\vecb}\to \R$ and that
\begin{equation*}
Y_s^{r,x,t,z,\vecb}=v_{\vecb}(s,X^{r,x}_s,t,z),\quad r\leq s\leq T.
\end{equation*}
Now, by Proposition~\ref{prop:VFbnd} the functions $v_{\vecb}$ are continuous and of polynomial growth. Finally, the verification theorem implies that the functions $v_{\vecb}$ are value functions for the stochastic control problem posed in Problem 1 in the sense that
\begin{align}\nonumber
v_{\vecb}(r,x,t,z)=&\sup_{\tau\in \mcT_t}\E\bigg[\int_r^{\tau\wedge T}\psi_{\vecb}\!\left(s,X_s^{r,x},R(z+(s-t)\vecb)\right)ds +\ett_{[\tau\geq T]}h_{\vecb}\!\left(X_T^{r,x},R(z+(T-t)\vecb)\right)\nonumber
\\
&+\ett_{[\tau < T]}\max_{\beta\in\mcJ^{-\vecb}}\left\{-c^{\vecb}_{\beta}(\tau) + v_{\beta}\!\left(\tau,X_\tau^{r,x},\tau,(z+(\tau-t)\vecb)^+\wedge \delta^{\beta}\right)\right\}\bigg].\label{ekv:VF}
\end{align}

\section{Limited feedback}
When searching for a numerical solution to Problem 1, by means of a lattice or a Monte Carlo approximation of the value function in \eqref{ekv:VF}, the curse of dimensionality will generally become apparent through an explosion in the computational burden as the number of units increase. To limit this effect we present an alternative, sub-optimal, scheme where only a part of the available state information is considered when making decisions.

Assume that, at time $t$, the system is operated in mode $\vecb\in\mcJ$ with $\zeta_t=\delta^{\vecb}$ when one or more units are intervened on giving us the new mode $\vecb'\in\mcJ^{-\vecb}$. The production in the period $[t,T]$ can then be written $p^{\vecb'}(\cdot,u)-\tilde R_{\vecb,\vecb'}(t,\cdot,u)$ where $u\in\mcU_t$ is the control applied in $[t,T]$,
\begin{align*}
p_i^{\vecb'}(s,u):=\:&\ett_{[b'_i=1]} \bigg\{\ett_{[t,\tau^i_1)}(s)\Cp_i+\sum_{j=1}^{\lfloor N_i/2 \rfloor}\ett_{[\tau^i_{2j},\tau^i_{2j+1})}(s)R_i\left((s-\tau^i_{2j})\wedge \delta_i\right)\bigg\}
\\
& +\ett_{[b'_i=0]} \sum_{j=1}^{\lceil N_i/2 \rceil}\ett_{[\tau^i_{2j-1},\tau^i_{2j})}(s)R_i\left((s-\tau^i_{2j-1})\wedge \delta_i\right)
\end{align*}
and
\begin{align*}
\left(\tilde R_{\vecb,\vecb'}(t,r,u)\right)_i=\ett_{[b'_i>b_i]}\ett_{[t,\tau^i_1)}(r)\{\Cp_i-R_{i}(r-t)\},
\end{align*}
for $i=1,\ldots,n$. The revenue without switching costs for the same period, for $X^{t,x}$, can then be written
\begin{align}\label{ekv:oc_tT}
\E\bigg[\int_t^T\psi_{\xi_s}\left(s,X^{t,x}_s,p^{\vecb'}(s,u)-\tilde R_{\vecb,\vecb'}(t,s,u)\right)ds + h_{\xi_T}\left(X^{t,x}_T,p^{\vecb'}(T,u)-\tilde R_{\vecb,\vecb'}(t,T,u))\right)\bigg].
\end{align}
If we define the \emph{delay revenue}
\begin{align*}
\Gamma^{\vecb}_{\vecb'}(t,x,u):=\E\bigg[&\int_t^T\{\psi_{\xi_r}(r,X_r^{t,x},p^{\vecb'}(r,u)-\tilde R_{\vecb,\vecb'}(t,r,u))-\psi_{\xi_r}(r,X_r^{t,x},p^{\vecb'}(r,u))\}dr
\\
&+h_{\xi_T}(X_T^{t,x},p^{\vecb'}(T,u)-\tilde R_{\vecb,\vecb'}(t,T,u))-h_{\xi_T}(X_T^{t,x},p^{\vecb'}(T,u))\bigg]
\end{align*}
then \eqref{ekv:oc_tT} can be written
\begin{align}\label{ekv:oc_tTnew}
\E\bigg[\int_t^T\psi_{\xi_r}\left(s,X^{t,x}_s,p^{\vecb'}(s,u)\right)ds + h_{\xi_T}\left(X^{t,x}_T,p^{\vecb'}(T,u))\right)+\Gamma^{\vecb}_{\vecb'}(t,x,u)\bigg].
\end{align}
This leads us to define a sequence of processes $((\hat Y^{t,x,\vecb,k}_s)_{0\leq s\leq T})_{k\geq 0}$ recursively as
\begin{align*}
\hat Y^{t,x,\vecb,0}_s:=\E\bigg[\int_s^{T}\psi_{\vecb}\left(r,X_r^{t,x},p^\vecb\right)dr+h_{\vecb}\left(X_T^{t,x},p^\vecb\right)\Big|\mcF_s\bigg],
\end{align*}
and
\begin{align*}
\hat Y^{t,x,\vecb,k}_s:=\esssup_{\tau\in\mcT_s} \E\bigg[&\int_s^{\tau\wedge T}\psi_{\vecb}(r,X_r^{t,x},p^\vecb)dr+\ett_{[\tau\geq T]}h_{\vecb}(X_T^{t,x},p^\vecb)
\\
&+\ett_{[\tau < T]}\max_{\beta\in \mcJ^{-\vecb}}\Big\{- c^{\vecb}_{\beta}(\tau)+\Gamma^{\vecb}_{\beta}(\tau,X^{t,x}_{\tau},u^{\diamond,k-1}_{\tau,X^{t,x}_{\tau},\beta})+\hat Y^{t,x,\beta,k-1}_{\tau}\Big\}\Big|\mcF_s\bigg],
\end{align*}
where $p^{\vecb}:=R(\delta^{\vecb})$ and the controls $u^{\diamond,k}_{t,x,\vecb}:=(\tau^{t,x,\vecb,k}_1,\ldots,\tau^{t,x,\vecb,k}_N;\beta^{t,x,\vecb,k}_1,\ldots,\beta^{t,x,\vecb,k}_N)\in\mcU_t^k$ are defined as follows. For each $\bbF$--stopping time $\tau$ we let
\begin{equation*}
D_\tau^{t,x,\vecb,k}:=\inf\Big\{s\geq \tau : \: \hat Y^{t,x,\vecb,k}_s=\max_{\beta\in \mcJ^{-\vecb}}\Big\{ -c^{\vecb}_{\beta}(s) +\Gamma^{\vecb}_{\beta}(s,X^{t,x}_{s},u^{\diamond,k-1}_{s,X^{t,x}_{s},\beta})+\hat Y^{t,x,\beta,k-1}_{s}\Big\}\Big\}.
\end{equation*}
We define the intervention times $\tau^{t,x,\vecb,k}_1,\ldots,\tau^{t,x,\vecb,k}_N$ and the corresponding sequence of active units $\beta^{t,x,\vecb,k}_1,\ldots,\beta^{t,x,\vecb,k}_N$ as $\tau^{t,x,\vecb,k}_1(=\tau_1):=D_t^{t,x,\vecb,k}$,
\begin{equation*}
\beta^{t,x,\vecb,k}_1(=\beta_1)\in\mathop{\arg\max}_{\beta \in \mcJ}\Big\{-c^{\vecb}_{\beta}(\tau_1) + \Gamma^{\vecb}_{\beta}(\tau_1,X^{t,x}_{\tau_1}, u^{\diamond,k-1}_{\tau_1,X^{t,x}_{\tau_1},\beta}) + \hat Y^{t,x,\beta,k-1}_{\tau_1}\Big\}
\end{equation*}
and continue with
\begin{equation*}
\tau^{t,x,\vecb,k}_j:=\tau^{\tau_1,X^{t,x}_{\tau_1},\beta_1,k-1}_{j-1} \quad \text{and} \quad \beta^{t,x,\vecb,k}_j:=\beta^{\tau_1,X^{t,x}_{\tau_1},\beta_1,k-1}_{j-1}
\end{equation*}
for $j=2,\ldots,N$.\\ 

For each $\vecb\in\mcJ$ and each $(t,x,z)\in [0,T]\times \R^m \times D_\zeta^{\vecb}$ we define the cost-to-go when applying the control $u\in\mcU_t$, as an extension of $J$ in the formulation of Problem~1,
\begin{align*}
J^{\vecb}(t,x,z;u):=\E\bigg[\int_t^T\psi_{\xi^{\vecb}_s}\left(s,X^{t,x}_s,R(\zeta^{t,z,\vecb}_s)\right)ds + h_{\xi^{\vecb}_T}\left(X^{t,x}_T,R(\zeta^{t,z,\vecb}_T)\right)-\sum_{j=1}^{N}c^{\beta_{j-1}}_{\beta_j}(\tau_j)\bigg].
\end{align*}

\begin{prop}
For all $\vecb\in\mcJ$, the sequence $((\hat Y^{t,x,\vecb,k}_s)_{0\leq s\leq T})_{k\geq 0}$ is an increasing sequence of \cadlag processes and satisfies
\begin{equation}
\lim_{k\to\infty} \hat Y^{t,x,\vecb,k}_s \leq Y_s^{t,x,0,\delta^{\vecb},\vecb},\label{ekv:hatYlimBND}
\end{equation}
$\Prob$--a.s.~for all $(s,t,x)\in [0,T]^2\times \R^m$. Furthermore, we have
\begin{align*}
\hat Y^{t,x,\vecb,k}_s=J^{\vecb}\left(s,X_s^{t,x},\delta^{\vecb};u^{\diamond,k}_{s,X_s^{t,x},\vecb}\right).
\end{align*}
\end{prop}

\noindent\emph{Proof.} Clearly, $(\hat Y^{t,x,\vecb,0}_s: 0\leq s\leq T)$ is \cadlag and $(\Gamma^{\vecb}_{\vecb'}(s,X^{t,x}_{s},\emptyset): 0\leq s\leq T)$ is continuous for each $\vecb,\vecb'\in \mcJ$, with $\vecb\neq \vecb'$. Hence, $\big(Y^{t,x,\vecb,1}_s+\int_0^s\psi_{\vecb}(r,X_r^{t,x},p^{\vecb})dr:0\leq s\leq T\big)$ is the Snell envelope of a \cadlag process and thus \cadlagp Now, assume that, for some $k\geq 1$, $(\hat Y^{t,x,\vecb,k}_s: 0\leq s\leq T)$ is \cadlag for all $\vecb\in \mcJ$. Then, $(\Gamma^{\vecb}_{\vecb'}(s,X^{t,x}_{s},u^{\diamond,k}_{s,X^{t,x}_{s},\vecb'}): 0\leq s\leq T)$ is \cadlag which implies that $(\hat Y^{t,x,\vecb,k}_s+\int_0^s\psi_{\vecb}(r,X_r^{t,x},p^{\vecb})dr: 0\leq s\leq T)$ is \cadlag as the Snell envelope of a \cadlag process.

From the proof of Theorem~1 in~\cite{BollanMSwitch1} and the Markov property we get
\begin{align*}
\hat Y^{t,x,\vecb,k}_s&=\esssup_{\tau\in\mcT_s} \E\bigg[\int_s^{\tau\wedge T}\psi_{\vecb}\left(r,X_r^{t,x},p^\vecb\right)dr+\ett_{[\tau \geq T]}h_{\vecb}\left(X_T^{t,x},p^\vecb\right)
\\
&\qquad+\ett_{[\tau < T]}\max_{\beta\in \mcJ^{-\vecb}}\Big\{- c^{\vecb}_{\beta}(\tau)+\Gamma^{\vecb}_{\beta}\left(\tau,X^{t,x}_{\tau},u^{\diamond,k-1}_{\tau,X^{t,x}_{\tau},\beta}\right)+\hat Y^{t,x,\beta,k-1}_{\tau}\Big\}\Big|\mcF_s\bigg]
\\
&=\E\bigg[\int_s^{\tau^{t,x,\vecb,k}_1\wedge T}\psi_{\vecb}\left(r,X_r^{t,x},p^\vecb\right)dr+\ett_{[\tau^{t,x,\vecb,k}_1\geq T]}h_{\vecb}\left(X_T^{t,x},p^\vecb\right)+\ett_{[\tau^{t,x,\vecb,k}_1 < T]}\Big\{- c^{\vecb}_{\beta^{t,x,\vecb,k}_1}(\tau^{t,x,\vecb,k}_1)
\\
&\qquad+\Gamma^{\vecb}_{\beta^{t,x,\vecb,k}_1}\left(\tau^{t,x,\vecb,k}_1,X^{t,x}_{\tau^{t,x,\vecb,k}_1}, u^{\diamond,k-1}_{\tau^{t,x,\vecb,k}_1,X^{t,x}_{\tau^{t,x,\vecb,k}_1},\beta^{t,x,\vecb,k}_1}\right)+\hat Y^{t,x,\beta^{t,x,\vecb,k}_1,k-1}_{\tau^{t,x,\vecb,k}_1}\Big\}\Big|\mcF_s\bigg]
\\
&=\ldots=
\\
&=\E\bigg[\int_s^{T}\sum_{j=0}^N \ett_{[\tau^{t,x,\vecb,k}_{j}\leq r < \tau^{t,x,\vecb,k}_{j+1}]}\psi_{\beta^{t,x,\vecb,k}_{j}}\left(r,X_r^{t,x},p^{\beta^{t,x,\vecb,k}_{j}}\right)dr + h_{\beta^{t,x,\vecb,k}_{N}}\left(X_T^{t,x},p^{\beta^{t,x,\vecb,k}_{N}}\right)
\\
&\qquad+ \sum_{j=1}^N\Big\{- c^{\beta^{t,x,\vecb,k}_{j-1}}_{\beta^{t,x,\vecb,k}_j}(\tau^{t,x,\vecb,k}_j) + \Gamma^{\beta^{t,x,\vecb,k}_{j-1}}_{\beta^{t,x,\vecb,k}_j} \Big(\tau^{t,x,\vecb,k}_j,X^{t,x}_{\tau^{t,x,\vecb,k}_j}, u^{\diamond,k-j}_{\tau^{t,x,\vecb,k}_j,X^{t,x}_{\tau^{t,x,\vecb,k}_j},\beta^{t,x,\vecb,k}_j}\Big)\Big\}\Big|\mcF_s\bigg]
\\
&=J^{\vecb}\left(s,X_s^{t,x},\delta^{\vecb};u^{\diamond,k}_{s,X_s^{t,x},\vecb}\right),
\end{align*}
with $(\tau^{t,x,\vecb,k}_{0},\beta^{t,x,\vecb,k}_{0})=(0,\vecb)$. Now, since $u^{\diamond,k}_{s,X_s^{t,x},\vecb}\in \mcU_s$ and $Y_s^{t,x,0,\delta^{\vecb},\vecb}=\esssup_{u\in\mcU_s}J^{\vecb}\left(s,X_s^{t,x},\delta^{\vecb};u\right)$, \eqref{ekv:hatYlimBND} follows.\qed\\


\subsection{Quadratic revenue\label{sec:quadREV}}
Consider now the problem of finding the control $u^*$ that maximizes \eqref{ekv:objFUN}
over all controls in $\mcU$ when, for each $b\in\mcJ$, $\psi_\vecb:[0,T]\times \R^k\times D_p\to \R$ and $h_\vecb$ are polynomials of degree two in $p$, \ie
\begin{align*}
\psi_{\vecb}(t,x,p)&=\psi_{\vecb}^0(t,x)+(\psi_{\vecb}^1(t,x))^\top p+p^\top \psi_{\vecb}^2(t,x) p,
\\
h_{\vecb}(x,p)&=h_{\vecb}^0(t,x)+(h_{\vecb}^1(x))^\top p+p^\top h_{\vecb}^2(x) p,
\end{align*}
where $a^\top$ is the transpose of the vector $a$ and, for all $\vecb\in\mcJ$, the functions $\psi_{\vecb}^0,h_{\vecb}^0:[0,T]\times \R^m\to \R$ and the components of $\psi_{\vecb}^1,h_{\vecb}^1:[0,T]\times \R^m\to \R^n$ and $\psi_{\vecb}^2,h_{\vecb}^2:[0,T]\times \R^m\to \R^{n\times n}$ are all locally Lipschitz continuous and of polynomial growth. Furthermore, we assume that the matrices $\psi_{\vecb}^2(t,x) $ and $h_{\vecb}^2(t,x)$ are both symmetric, for all $(t,x)\in [0,T]\times \R^m$.
The delay revenue can then be written
\begin{align*}
\Gamma^{\vecb}_{\vecb'}(t,x,u)=\E\bigg[&\int_t^T\Big\{-\left((\psi_{\xi^{\vecb'}_r}^1(r,X_r^{t,x}))^\top +2(p^{\vecb'}(r,u))^\top \psi_{\xi^{\vecb'}_r}^2(r,X_r^{t,x})\right)\tilde R_{\vecb,\vecb'}(t,r,u)
\\
&+(\tilde R_{\vecb,\vecb'}(t,r,u))^\top \psi^2_{\xi^{\vecb'}_r}(r,X_r^{t,x})\tilde R_{\vecb,\vecb'}(t,r,u)\Big\}dr
\\
&-\left(h^1_{\xi^{\vecb'}_T}(X_T^{t,x})+2(p^{\vecb'}(T,u))^\top h^2_{\xi^{\vecb'}_T}(X_T^{t,x})\right)\tilde R_{\vecb,\vecb'}(t,T,u)
\\
&+(\tilde R_{\vecb,\vecb'}(t,T,u))^\top h^2_{\xi^{\vecb'}_T}(X_T^{t,x})\tilde R_{\vecb,\vecb'}(t,T,u)\bigg].
\end{align*}
Exploiting the simple structure of this formulation, we will in what follows show that efficient numerical algorithms can be built to approximate the expected revenue and the corresponding control $u^\diamond$.

\subsection{Numerical solution scheme}
In this section we present a numerical scheme that approximates $J^{\vecb}\left(t,x,\delta^{\vecb};u^{\diamond,k}_{t,x,\vecb}\right)$ when the $\psi_{\vecb}$ and the $h_{\vecb}$ are quadratic polynomials in $p$.

We start by going from continuous to discrete time by introducing the grid $\Pi=\{t_0,t_1,\ldots,t_{N_{\Pi}}\}$, with $t_l= l\Delta t$ for $l=0,\ldots, N_{\Pi}$, where $\Delta t=T/N_{\Pi}$. To get a discrete time problem we apply a Bermudan options approximation and reduce the set of stopping times in the admissible controls by restricting interventions to grid points, \ie for all discretized intervention times $\bar\tau_j$ we have $\bar\tau_j\in \Pi$.

Let $\hat u^{\diamond}_{t,x,\vecb}:=(\bar\tau^{t,x,\vecb}_1,\ldots,\bar\tau^{t,x,\vecb}_N;\beta^{t,x,\vecb}_1,\ldots,\beta^{t,x,\vecb}_N)$, be the discreet-time version of the limited feedback control proposed above and let $(\hat\xi^{t,x,\vecb}_{s}:s\in \Pi\cap [t,T])$ be the corresponding evolution of the operating mode. We define the discrete time value function
\begin{align*}
\hat v^{\Pi}_\vecb(t_l,x):=\E\bigg[&\sum_{k=l}^{N_\Pi-1} \psi_{\hat\xi^{t_l,x,\vecb}_{t_k}}\left(t_k,X_{t_k}^{t_l,x},p^\vecb(t_k,\hat u^{\diamond}_{t_l,x,\vecb})\right)\Delta t
\\&+h_{\hat\xi^{t_l,x,\vecb}_{T}}\left(X_{T}^{t_l,x},p^\vecb(T,\hat u^{\diamond}_{t_l,x,\vecb})\right) -\sum_{j=1}^N c^{\beta^{t_l,x,\vecb}_{j-1}}_{\beta^{t_l,x,\vecb}_j}(\bar\tau^{t_l,x,\vecb}_j)\bigg],
\end{align*}
with $(\bar\tau^{t_l,x,\vecb}_{0},\beta^{t_l,x,\vecb}_0)=(0,\vecb)$. Then the functions $\hat v^{\Pi}_\vecb:\Pi\times \R^m\to\R$ satisfy the recursion
\begin{align}
\hat v^{\Pi}_\vecb (T,x)&=h_{\vecb}(x,p^{\vecb}), \label{ekv:diskBellmanT}\\
\hat v^{\Pi}_\vecb(t_l,x)&=\max \limits_{\beta\in \mcJ}\left\{\psi_{\beta}(t_l,x,p^{\vecb}\wedge p^{\beta})\Delta t-c^{\vecb}_{\beta}(t_{l})+\E\left[\hat\Gamma^{\vecb}_{\beta}(t_{l+1},X^{t_l,x}_{t_{l+1}})+\hat v^{\Pi}_\beta(t_{l+1},X^{t_l,x}_{t_{l+1}}) \right]\right\}, \label{ekv:diskBellman}
\end{align}
where, for each $\vecb,\vecb'\in\mcJ$, the discrete-time delay revenue $\hat\Gamma^{\vecb}_{\vecb'}:\Pi\times \R^m \to \R$ is given by
\begin{align*}
\hat\Gamma^{\vecb}_{\vecb'}(t_l,x)=\sum_{k={l}}^{N_\Pi}\Big\{&-\sum_{i\in\mcI(\vecb')\setminus\mcI(\vecb)} \Lambda^1_{\vecb',i}(t_{l-1},x;t_k)\tilde R_i(t_k-t_{l-1})
\\
&+\sum_{i,j\in\mcI(\vecb')\setminus\mcI(\vecb)}\Lambda^2_{\vecb',i,j}(t_{l-1},x;t_k)\tilde R_i(t_k-t_{l-1})\tilde R_j(t_k-t_{l-1})\Big\}\Delta t
\end{align*}
with, for all $i\in\mcI(\vecb')$,
\begin{align}
\Lambda^1_{\vecb',i}(t_l,x;t_k)&:=\E\Big[\ett_{[\bar \tau_1^{i,t_l,x,\vecb'}>t_k]}\big\{(\psi_{\hat\xi^{t_l,x,\vecb'}_{t_k}}^1(t_k,X_{t_k}^{t_l,x}))_i +2((p^{\vecb'}(t_k,\hat u^{\diamond}_{t_l,x,\vecb'}))^\top \psi^2_{\hat\xi^{t_l,x,\vecb'}_{t_k}}(t_k,X_{t_k}^{t_l,x}))_i\big\}\Big],\label{ekv:Lamb1}
\\
\Lambda^1_{\vecb',i}(t_l,x;T)&:=\E\Big[\ett_{[\bar \tau_1^{i,t_l,x,\vecb'}>T]}\big\{(h_{\hat\xi^{t_l,x,\vecb'}_{T}}^1(X_{T}^{t_l,x}))_i +2((p^{\vecb'}(T,\hat u^{\diamond}_{t_l,x,\vecb'}))^\top h^2_{\hat\xi^{t_l,x,\vecb'}_{T}}(X_{T}^{t_l,x}))_i\big\}\Big]\label{ekv:Lamb1T}
\end{align}
and, for all $i,j\in\mcI(\vecb')$,
\begin{align}
\Lambda^2_{\vecb',i,j}(t_l,x;t_k)&:=\E\Big[\ett_{[\bar \tau_1^{i,t_l,x,\vecb'}>t_k]} \ett_{[\bar \tau_1^{j,t_l,x,\vecb'}>t_k]} (\psi^2_{\hat\xi^{t_l,x,\vecb'}_{t_k}}(t_k,X_{t_k}^{t_l,x}))_{i,j}\big\}\Big],\label{ekv:Lamb2}
\\
\Lambda^2_{\vecb',i,j}(t_l,x;T)&:=\E\Big[\ett_{[\bar \tau_1^{i,t_l,x,\vecb'}>T]} \ett_{[\bar \tau_1^{j,t_l,x,\vecb'}>T]} (h^2_{\hat\xi^{t_l,x,\vecb'}_{T}}(X_T^{t,x}))_{i,j}\big\}\Big].\label{ekv:Lamb2T}
\end{align}
\begin{rem}
In each $t_l\in\Pi\setminus \{T\}$ the recursion \eqref{ekv:diskBellman} evaluates the optimal action to take at the present time, take the action and move to the next time $t_{l+1}\in\Pi$ where the process is repeated. Having arrived at the conclusion that $\beta\in\mcJ$ is the optimal action, we know that the present production is $p^{\vecb}\wedge p^\beta$ as turnoffs are immediate while increasing the output requires ramping. This is why $\psi$ is evaluated in $p^{\vecb}\wedge p^\beta$.
\end{rem}

\subsubsection{Recursions for $\Lambda^1_{\vecb,i}$ and $\Lambda^2_{\vecb,i,j}$}
At each time step, starting at $t_l=T$ and moving backwards, we obtain the expected revenue to-go by solving \eqref{ekv:diskBellmanT} and \eqref{ekv:diskBellman}. This gives us, for each $\vecb\in\mcJ$ and each $t_l\in\Pi\setminus\{T\}$, the $\mcF_{t_l}$--measurable optimal actions as a selection of
\begin{align}
\beta^*(\vecb,t_l,x)\in\mathop{\arg\max} \limits_{\beta\in \mcJ}\left\{\psi_{\beta}(t_l,x,p^{\vecb}\wedge p^{\beta})\Delta t-c^{\vecb}_{\beta}(t_{l})+\E\left[\hat\Gamma^{\vecb}_{\beta}(t_{l+1},X^{t_l,x}_{t_{l+1}})+\hat v^{\Pi}_\beta(t_{l+1},X^{t_l,x}_{t_{l+1}}) \right]\right\}, \label{ekv:bstar}
\end{align}
from which we deduce that the intervention times satisfy $\{\bar \tau_1^{i,t_l,x,\vecb}>t_l\}=\{(\beta^*(\vecb,t_l,x))_i= b_i\}$, for $i=1,\ldots,n$. As we will see, knowledge of whether these events occur is enough to compute $\Lambda^1_{\vecb,i}$ and $\Lambda^2_{\vecb,i,j}$ in a recursive manner.

Let us start with the simpler $\Lambda^2_{\vecb,i,j}$, where $\vecb\in\mcJ$ and $i,j\in\mcI(\vecb)$. First, if any of the events $\{\bar \tau_1^{i,t_l,x,\vecb}=t_l\}$ and $\{\bar \tau_1^{j,t_l,x,\vecb}=t_l\}$ occur, then \eqref{ekv:Lamb2} and \eqref{ekv:Lamb2T} immediately give $\Lambda^2_{\vecb,i,j}(t_l,x;t_k)=0$, for all $t_k\in\Pi$ with $t_k\geq t_l$.

Assume instead that $\bar \tau_1^{i,t_l,x,\vecb}>t_l$ and $\bar \tau_1^{j,t_l,x,\vecb}>t_l$, $\Prob$-a.s. Then
\begin{equation*}
\Lambda^2_{\vecb,i,j}(t_l,x;t_l)=\left\{\begin{array}{ll}
(\psi^2_{\beta^*}(t_l,x))_{i,j}, & \text{for }t_l\in \Pi\setminus \{T\}, \\
(h^2_{\vecb}(x))_{i,j}, & \text{for }t_l=T,
\end{array}\right.
\end{equation*}
and
\begin{equation*}
\Lambda^2_{\vecb,i,j}(t_l,x;t_k)=\E\big[\Lambda^2_{\beta^*,i,j}(t_{l+1},X_{t_{l+1}}^{t_l,x};t_k)\big],\quad\text{for }t_k>t_l.
\end{equation*}

For $\Lambda^1_{\vecb,i}$ the situation is just slightly more involved, as these depend on future values of the optimal output vector $p^{\vecb}(t_k,\bar u^{\diamond}_{t_l,x,\vecb})$.
As above we note that, whenever $\bar \tau_1^{i,t_l,x,\vecb}=t_l$, equations \eqref{ekv:Lamb1} and \eqref{ekv:Lamb1T} give $\Lambda^1_{\vecb,i}(t_l,x;t_k)=0$, for all $t_k\geq t_l$.

Let us thus assume that $\bar \tau_1^{i,t_l,x,\vecb}>t_l$. From \eqref{ekv:diskBellmanT} and \eqref{ekv:diskBellman} we note that we must have
\begin{equation*}
\Lambda^1_{\vecb,i}(t_l,x;t_l)=\left\{\begin{array}{ll}
(\psi_{\beta^*}^1(t_k,x) + 2(p^{\vecb}\wedge p^{\beta^*})^\top \psi^2_{\beta^*}(t_k,x))_i, & \text{for }t_l\in \Pi\setminus \{T\}, \\
(h_{\vecb}^1(x) + 2(p^{\vecb})^\top h^2_{\vecb}(x))_i, & \text{for }t_l=T.
\end{array}\right.
\end{equation*}
In the recursion for $\Lambda^1_{\vecb,i}(t_l,x;t_l)$ we have to consider the fact that $p^{\beta^*}(t_k,\hat u^{\diamond}_{t_l,x,\beta^*})$ depends on future control actions. In particular $(p^{\vecb}(t_k,\hat u^{\diamond}_{t_l,x,\vecb}))_j=(p^{\beta^*}(t_k,\hat u^{\diamond}_{t_l,x,\beta^*}))_j-\ett_{[\bar\tau_1^{j,t_l,x,\beta^*}>t_k]}\tilde R_j(t_k-t_l)$, for all $j\in \mcI(\beta^*)\setminus \mcI(\vecb)$. Hence,
\begin{equation*}
\Lambda^1_{\vecb,i}(t_l,x;t_k)=\E\big[\Lambda^1_{\beta^*,i}(t_{l+1},X_{t_{l+1}}^{t_l,x};t_k)\big]-\sum_{j\in \mcI(\beta^*)\setminus \mcI(\vecb)}\Lambda^2_{\beta^*,i,j}(t_l,x;t_k)\tilde R_j(t_k-t_l),\quad\text{for }t_k>t_l.
\end{equation*}


\section{Numerical example}
\begin{figure}[h!]
  \centering
  \includegraphics[width=0.75\textwidth]{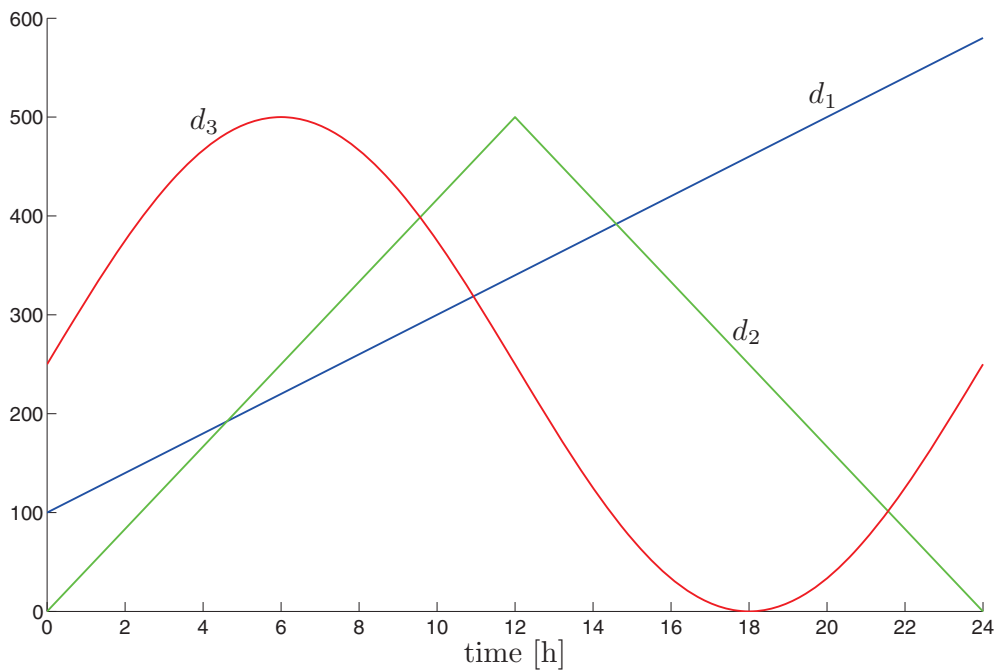}
  \caption{The three different $d$ that are investigated.}\label{fig:mFUNs}	
\end{figure}

In the numerical example we will consider a tracking problem where an operator wants to minimize
\begin{equation}\label{ekv:TRprob}
J(u):=\E\bigg[\int_0^T (f_{\textit{pen}}(X_t-\sum_{i=1}^n p_i(t))^2+(c^f)^{\!\top} p(t)) dt+f_{\textit{pen},T}(X_T-\sum_{i=1}^n p_i(T))^2+\sum_{j=1}^N c^{\beta_{j-1}}_{\beta_j}(\tau_j)\bigg]
\end{equation}
over a period of $T=24$ hours, where $f_{\textit{pen}},f_{\textit{pen},T}>0$ are penalization coefficients, $c^f\in\R_+^n$ is the marginal production cost in the different units and the switching costs are constant. The signal $(X_t:0\leq t\leq T)$ to be tracked is given by the sum, $X_t=d(t)+Z_t$ of a deterministic forecast $(d(t):0\leq t\leq T)$ and an Ornstein-Uhlenbeck process that solves the SDE
\begin{align*}
dZ_t&=-a Z_t dt + \sigma dW_t,\quad \text{for } t\in[0,T]
\\
Z_0&=x_0-d(0),
\end{align*}
where $a=0.01$ and $\sigma=10$. We will investigate the performance of the limited feedback control $u^\diamond$ for three different shapes of the forecast $d(t)$,
\begin{align*}
d_1(t)&:=100+20t,\\
d_2(t)&:=500\left(1-\tfrac{2}{T}|t-T/2|\right),\\
d_3(t)&:=250(1+\sin(2\pi t/T)),
\end{align*}
that are depicted in Fig.~\ref{fig:mFUNs}.

We assume that the operator has at her disposal a set of six production units whose data is summarized in Table~\ref{tab:unitDATA}, where $c^f_i$ is the marginal production cost in Unit $i$ and the associated ramp function is defined through the constants $0\leq \delta'_i<\delta_i$ and $\Cp_i$ as
\begin{equation*}
R_i(s)=1_{[\delta'_i,\delta_i]}(s)\frac{s-\delta'_i}{\delta_i-\delta'_i}\Cp_i,\quad \text{for }i=1,\ldots,6.
\end{equation*}

\begin{table}[h!]\centering \setlength{\tabcolsep}{10pt}
\begin{tabular}{| c | c | c | c | c | c | c |}\hline
  $i$ & $\Cp_i$ & $c_i^{\on}$ & $c_i^{\off}$ & $c_i^f$ & $\delta'_i$ & $\delta_i$ \\
  \hline
   1  &  150  &  5000 &  2000   &    3   &    2    &    7 \\ \hline
   2  &  125  &  3000 &  2000   &    4   &    1    &    6 \\ \hline
   3  &  100  &  2000 &  1500   &    4   &    2    &    5 \\ \hline
   4  &   75  &  1500 &  1000   &    4   &    1    &    4 \\ \hline
   5  &   50  &   750 &  1000   &    5   &    1    &    3 \\ \hline
   6  &   25  &   500 &  1000   &    7   &    1    &    2 \\ \hline
\end{tabular}
\caption{Data for the production units in the example.}\label{tab:unitDATA}
\end{table}

Equation \eqref{ekv:TRprob} can be written
\begin{equation*}
J(u)=\E\bigg[\int_0^T \Big([X_t\:\: (p(t))^{\!\top}]Q\left[\begin{array}{c}X_t \\ p(t)\end{array}\right]+(c^f)^{\!\top} p(t)\Big)dt+[X_T\: \: (p(T))^{\!\top}]M\left[\begin{array}{c}X_T \\ p(T)\end{array}\right]-\sum_{j=1}^N c^{\beta_{j-1}}_{\beta_j}(\tau_j)\bigg],
\end{equation*}
where $Q$ and $M$ are symmetric matrices. Hence, the problem of finding an efficient control scheme fits in the quadratic setting described in Section~\ref{sec:quadREV}.

We solve the problem for constants $f_{\textit{pen}}=0.1$ and $f_{\textit{pen},T}=0.3$ and using three different sets of available units $F_1:=\{3,5\}$, $F_2:=\{2,4,6\}$ and $F_3:=\{1,2,3,4,5,6\}$ for each of the three different forecasts.

The problem is numerically solved by means of a Markov-Chain approximation of the process $(X_t:0\leq t\leq T)$ as prescribed in~\cite{numSCbok}. We use a time-discretization with $N_{\Pi}=241$ points and discretize the state space of $(X_t:0\leq t\leq T)$ using 201 grid-points.

With this dicsretization, the numerical solution was obtained in 4, 18 and 720 seconds for the limited feedback algorithm. For the fully augmented solution method the first two settings with two and three units where solved in around 220 and 12000 seconds, respectively (it seemed computationally impossible to obtain a solution with the full system of six units).

\begin{figure}[h!]
  \centering
  \includegraphics[width=0.75\textwidth]{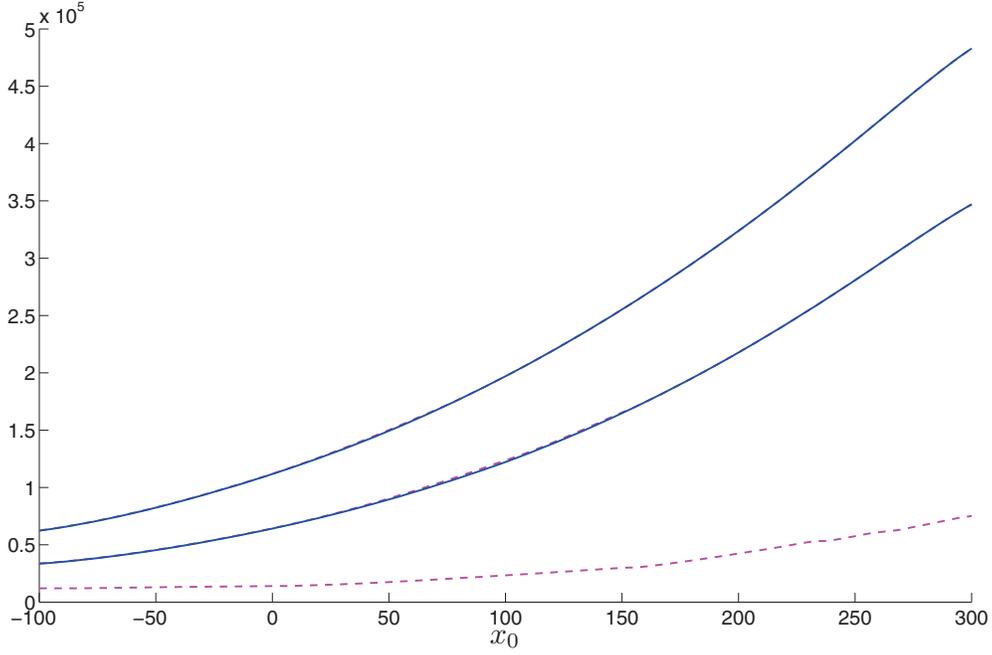}
  \caption{Cost-to-go at $t=0$ for limited feedback algorithm (solid blue) and minimal operation cost (dashed magenta) with forecast $d_1$.}\label{fig:VFun1}
\end{figure}

\begin{figure}[h!]
  \centering
  \includegraphics[width=0.75\textwidth]{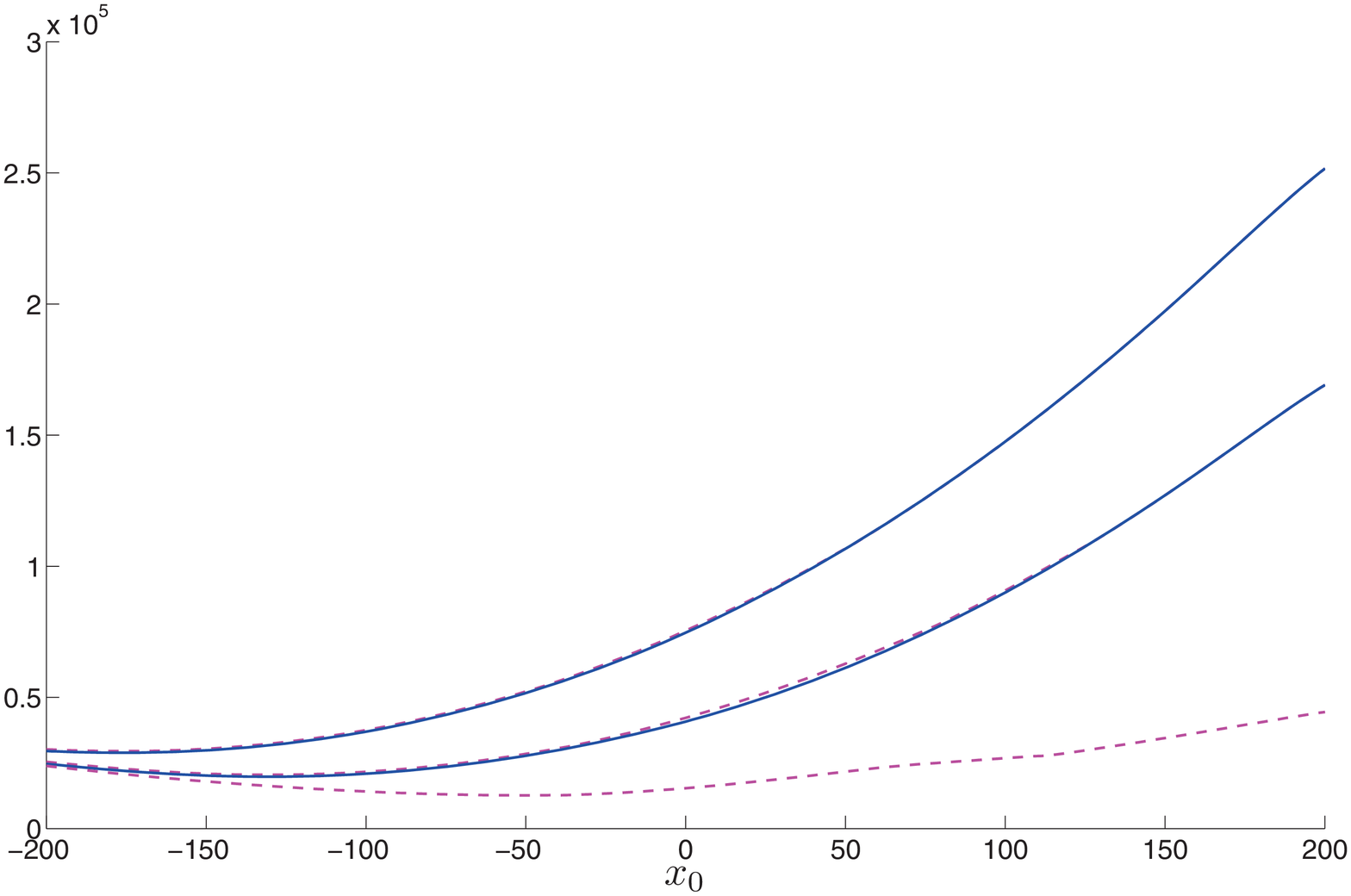}
  \caption{Cost-to-go at $t=0$ for limited feedback algorithm (solid blue) and minimal operation cost (dashed red).}\label{fig:VFun2}
\end{figure}

\begin{figure}[h!]
  \centering
  \includegraphics[width=0.75\textwidth]{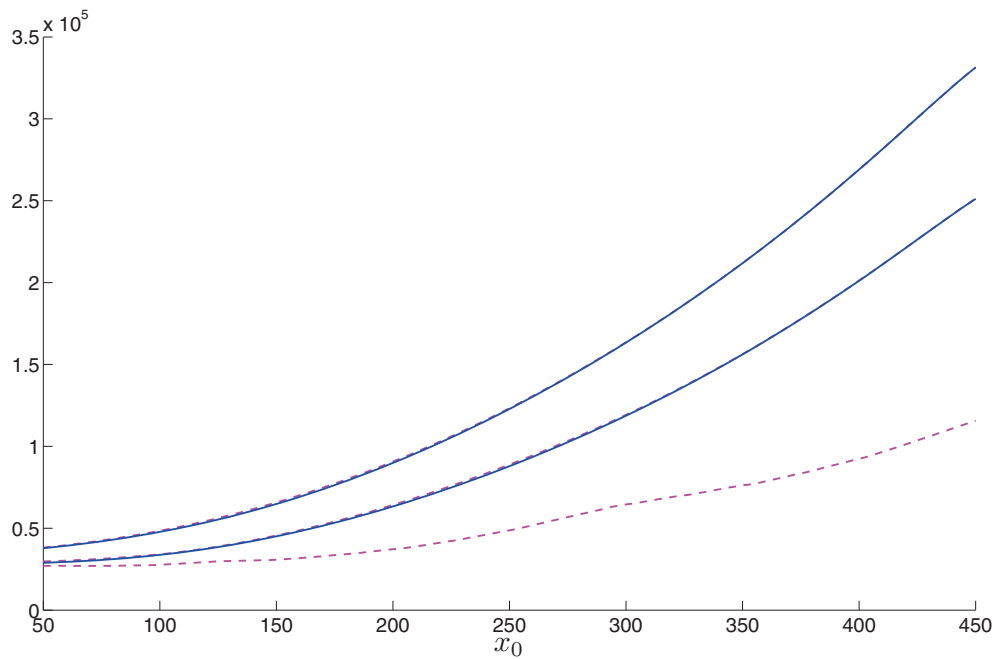}
  \caption{Cost-to-go at $t=0$ for limited feedback algorithm (solid blue) and minimal operation cost (dashed red).}\label{fig:VFun3}
\end{figure}

Figures~\ref{fig:VFun1}-\ref{fig:VFun3} show the expected operation costs at time zero for the limited feedback approach (solid blue lines) and the corresponding minimal operation costs obtained by state space augmentation (dashed magenta lines), for the three different forecasts. In all cases the expected operation costs decreased with more units, in particular the expected operation cost with units $\{3,5\}$ was always higher than the expected operation cost with units $\{2,4,6\}$.

\begin{figure}[h!]
  \centering
  \includegraphics[width=0.75\textwidth]{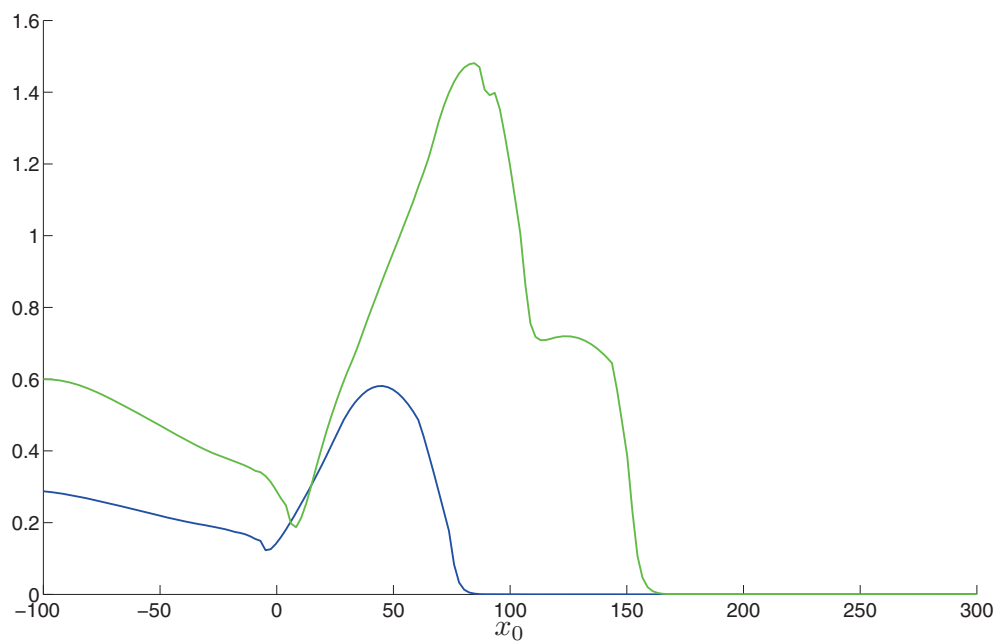}
  \caption{Relative errors (in \%) with $d_1$.}\label{fig:relERR1}
\end{figure}

\begin{figure}[h!]
  \centering
  \includegraphics[width=0.75\textwidth]{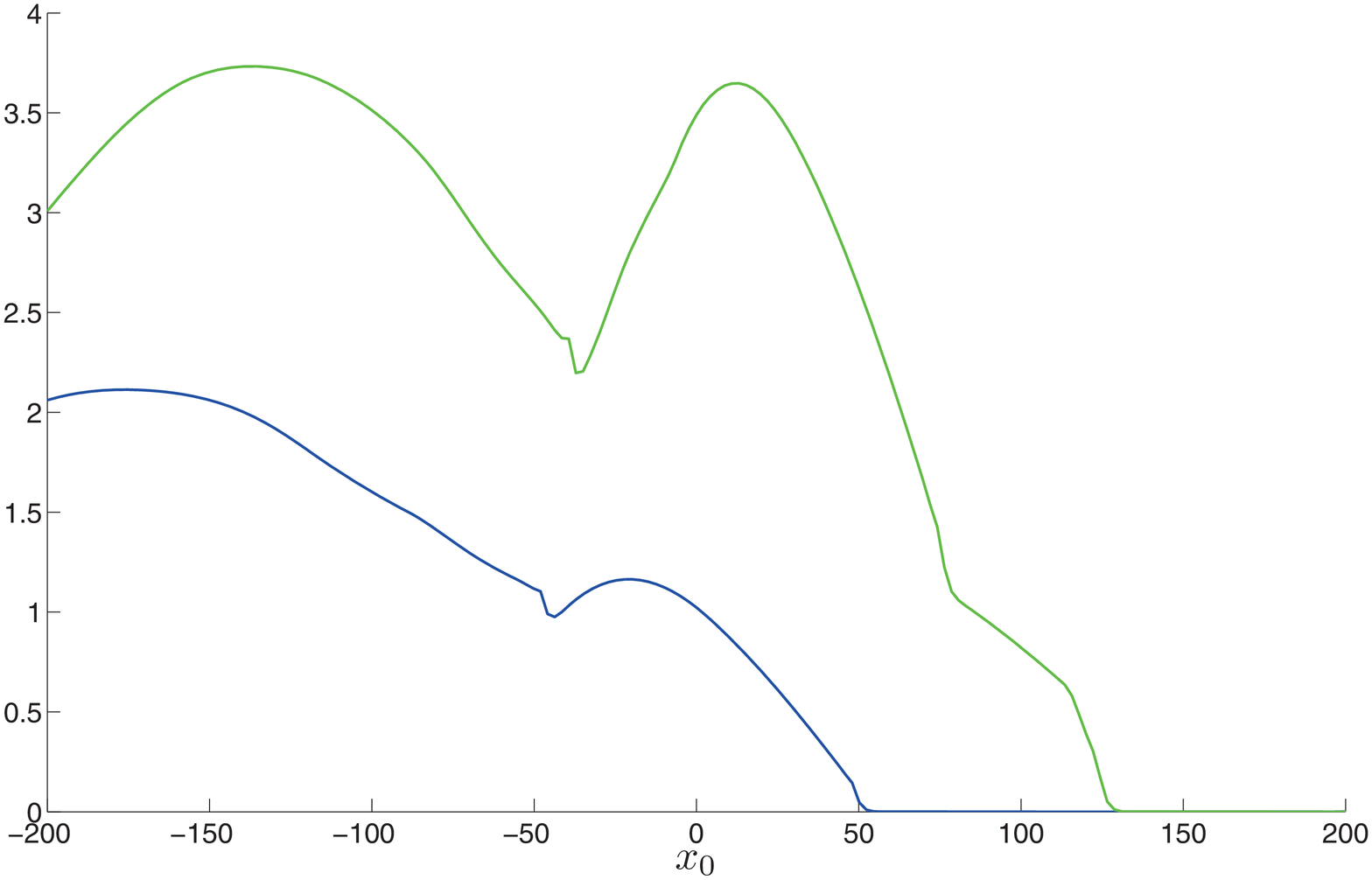}
  \caption{Relative errors (in \%) with $d_2$.}\label{fig:relERR2}
\end{figure}

\begin{figure}[h!]
  \centering
  \includegraphics[width=0.75\textwidth]{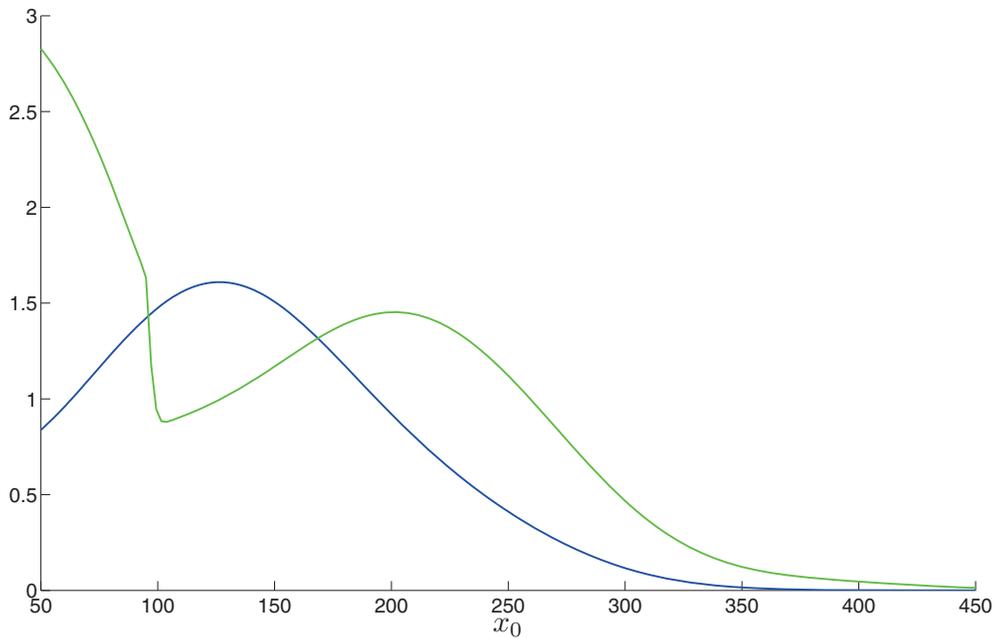}
  \caption{Relative errors (in \%) with $d_3$.}\label{fig:relERR3}
\end{figure}

In Figures~\ref{fig:relERR1}-\ref{fig:relERR3} the relative error of the limited feedback approximation is plotted for the three different forecasts. Here, we define the relative error as the function
\begin{align*}
e_{rel}(x):=100\left(\frac{\hat v_\bold{0}^{\Pi}(0,x)}{v_{\bold{0}}^{\Pi}(0,x,0,0)}-1\right),
\end{align*}
where $v_{\bold{0}}^{\Pi}$ is the discretized version of $v_{\bold{0}}$. In the figures the blue lines are the relative errors with units $\{3,5\}$ and the green lines are the relative errors with units $\{2,4,6\}$.

Note that the level of sub-optimality induced by the limited feedback approximation depends on the properties of the process $(X_t: 0\leq t\leq T)$ but also on the available units. The seemingly higher error with three units ($F_2$) compared to with two units ($F_1$) can, however, be partially explained by the lower operation cost for $F_2$ leading to a higher weight of the absolute error in the relative error.

\bibliographystyle{plain}
\bibliography{LimFEED_ref}
\end{document}